\documentclass[12pt]{amsart}
\usepackage{amscd,amssymb}
\newcommand{\Aut}{\operatorname{Aut}}

\newcommand{\bx}{\mathbf{x}}
\newcommand{\cft}{CohFT}
\newcommand{\cfts}{CohFTs}
\newcommand{\CP}[1]{\mathbb{CP}^{#1}}
\newcommand{\End}{\mathcal{E}nd}
\newcommand{\G}{\mathcal{G}} 	

\newcommand{\HH}{\mathcal{H}}	

\newcommand{\Hp}{\operatorname{\tilde{H}}}

\newcommand{\LL}{\mathcal{L}}	
\newcommand{\mf}{\mathbf{m}}    
\newcommand{\M}{\overline{\MM}}	
\newcommand{\MM}{\mathcal{M}}	

\newcommand{\fprod}{\widehat{\prod}} 

\newcommand{\pf}{\mathbf{p}}

\newcommand{\qh}{\widetilde{q}}
\newcommand{\R}{\mathcal{R}}	
\newcommand{\SC}{\mathcal{S}}
\renewcommand{\sf}{\mathbf{s}}
\newcommand{\st}{\operatorname{st}}
\newcommand{\tf}{\mathbf{t}}

\newcommand{\Deltat}{\widetilde{\Delta}}
\newcommand{\dd}{\partial} 
\newcommand{\delf}{\boldsymbol{\delta}}

\newcommand{\kaf}{\boldsymbol{\kappa}} 
\newcommand{\kao}{\widehat{\kappa}} 
\newcommand{\pit}{\widetilde{\pi}} 
\newcommand{\psio}{\widehat{\psi}} 
\newcommand{\tauf}{\boldsymbol{\tau}}
\newcommand{\zef}{\boldsymbol{0}} 

\newcommand{\nc}{{\mathbb{C}}} 	
\newcommand{\nz}{{\mathbb{Z}}}	

\newcommand{\la}{\langle}  
\newcommand{\ra}{\rangle}  

\def\lift#1#2{
  \dimen0 = \unitlength 
  \multiply\dimen0 by #1 \divide \dimen0 by 2
  \dimen1 = \dimen0 
  \multiply \dimen1 by 7 \divide \dimen1 by 10
  \raise\dimen1
     \hbox{\hskip 0.3cm ${\vbox to \dimen0{}}$ \enspace #2}}

\newtheorem{thm}{Theorem}[section]
\newtheorem{lm}[thm]{Lemma}
\newtheorem{prop}[thm]{Proposition}
\newtheorem{crl}[thm]{Corollary}

\theoremstyle{definition}

\newtheorem{df}[thm]{Definition}

\theoremstyle{remark}

\newtheorem{nota}{Notation} 
\newtheorem{rem}{Remark}  
\newtheorem{ack}{Acknowledgment}

\begin{document}
\addtocounter{section}{-1}

\title[Intersection Numbers]
{Intersection Numbers on the Moduli Spaces of Stable Maps in Genus $0$} 

\author
[A. Kabanov]{Alexandre Kabanov}
\address
{Max Planck Institut f\"ur Mathematik, Gottfried Claren Str. 26, 53225
Bonn, GERMANY and Department of Mathematics, Michigan State
University, Wells Hall,
East Lansing, MI 48824-1027, USA}
\email{kabanov@math.msu.edu}

\author
[T. Kimura]{Takashi Kimura}
\address
{Max Planck Institut f\"ur Mathematik, Gottfried Claren Str. 26, 53225
Bonn, GERMANY and Department of Mathematics, 111 Cummington Street, Boston
University, Boston, MA 02215, USA} 
\email{kimura@math.bu.edu}

\date{\today}

\begin{abstract} 
Let $V$ be a smooth, projective, convex variety. We define
tautological $\psi$ and $\kappa$ classes on the moduli space of stable
maps $\M_{0,n}(V)$, give a (graphical) presentation for these classes
in terms of boundary strata, derive differential equations for the
generating functions of the Gromov-Witten invariants of $V$ twisted by
these tautological classes, and prove that these intersection numbers
are completely determined by the Gromov-Witten invariants of $V$. This
results in families of Frobenius manifold structures on the cohomology
ring of $V$ which includes the quantum cohomology as a special case.
\end{abstract}

\maketitle


\section{Introduction}
\label{intro}

There has recently been a great deal of interest in $\M_{g,n}(V)$, the moduli
space of stable maps of genus $g$ with $n$ marked points into a smooth,
projective variety $V$, an object whose construction was envisioned by
Kontsevich as the proper algebro-geometric setting for Gromov-Witten
invariants and quantum cohomology. The space $\M_{g,n}(V)$ is a 
compactification of the moduli space of holomorphic maps 
from a genus $g$ Riemann surface with $n$ marked points to $V$ by allowing
the surfaces to degenerate forming double points away from the marked points
provided that a stability condition is satisfied. The moduli space of stable
maps $\M_{g,n}(V)$ should be regarded as a generalization of the moduli space
of stable curves $\M_{g,n}$, due to Deligne-Mumford, as it reduces to the
latter in the special case where $V$ is a point.  However, unlike the moduli
space of curves, $\M_{g,n}(V)$ is generally a singular space (a
Deligne-Mumford stack) when the relevant obstruction bundle fails to
vanish. Intersection theory on stacks requires the usage of the virtual
fundamental class \cite{LT,BF, BMa} in  lieu of the topological one. However,
when $V$ is a smooth, projective, convex variety $\M_{0,n}(V)$ are complex
orbifolds \cite{FP}. We shall henceforth assume $V$ satisfies these
conditions and we restrict ourselves to genus zero.

The spaces $\M_{0,n}(V)$ are equipped with evaluation morphisms
$\mathrm{ev}_i\,:\,\M_{0,n}(V)\,\to\,V$ for $i\,=\,1,\ldots,n$ which
evaluates the stable map on the $i^\mathrm{th}$ puncture. The Gromov-Witten
invariants of $V$ are the intersection numbers of the pull back of cohomology
classes on $V$ via these evaluation morphisms. The Gromov-Witten invariants
often solve enumerative problems in algebraic geometry, {\sl e.g.}\ the
number of rational curves in $\CP{2}$ counted with suitable multiplicity
\cite{KM1} (see also \cite{DFI,FP}). The Gromov-Witten invariants of $V$
endow $H^\bullet(V)$ with a deformed cup product resulting in the structure
of quantum cohomology (see also \cite{RT}) which makes $H^\bullet(V)$ into a
cohomological field theory (\cft) \cite{KM1,Ma} (or 
equivalently, a (formal) Frobenius manifold \cite{Du,Hi,Ma}). The
Gromov-Witten invariants can be encoded in a generating function (called the
potential) on $H^\bullet(V)$ which must obey the WDVV
(Witten-Dijkgraaf-Verlinde-Verlinde) equations. The WDVV equations allows the
recursive computation of the Gromov-Witten invariants of certain homogeneous
spaces, \emph{e.g.}\ projective spaces, Grassmannians, etc. \cite{KM1,FP,DFI}

Another set of cohomology classes on $\M_{0,n}(V)$ are those associated to
its universal curve. The tautological $\psi$ classes are the first Chern
classes of tautological line bundles on $\M_{0,n}(V)$. Certain intersection
numbers of the $\psi$ classes and the pull back of cohomology classes on $V$
via the evaluation morphisms are the called Gromov-Witten invariants twisted
by the $\psi$ classes.  These twisted invariants are of great interest
as the associated potential is conjectured to obey a highest weight condition
for the Virasoro Lie algebra \cite{E,E2}, a conjecture which has some
intriguing evidence behind it \cite{Ba,Pa2}. This conjecture has been proven
in the case where $V$ is a point and is a consequence of Kontsevich's proof
of the Witten conjecture \cite{Lo,Ko,Di}. There, the highest weight condition
gives rise to recursion relations which completely determine these
intersection numbers in all genera. 

In this paper, we begin by exploiting the canonical stratification of
$\M_{0,n}(V)$ by complex orbifolds which generalize the stratification of
$\M_{0,n}$. These strata are indexed by decorated trees which can be used to
develop a graphical language describing cohomology classes on
$\M_{0,n}(V)$. This graphical language is particularly useful in describing
the behavior of these classes under restriction to the strata and under the
push forward and pull back with respect to the forgetful map $\M_{0,n+1}(V)
\, \to\,\M_{0,n}(V)$.

We then introduce the tautological $\psi$ and $\kappa$ classes on
$\M_{0,n}(V)$ generalizing those classes from $\M_{0,n}$. The definition of
the $\kappa$ classes is new and has the novel feature that their proper 
definition involves the pull backs via the evaluation morphisms of the
cohomology classes on $V$. We prove restriction properties for these classes
to the boundary strata and derive a (graphical) presentation of the $\psi$
and $\kappa$ classes in terms of boundary strata generalizing the
presentation for the same classes in \cite{KK} on the moduli space of curves.

We prove that the Gromov-Witten classes twisted by the $\psi$ and $\kappa$
classes endows $H^\bullet(V)$ with a family of \cft\ structures by
pushing them forward to the moduli space of curves. These \cft\  structures
are deformations of the cup product containing, as a special case, the
product in quantum cohomology.

Finally, we prove recursion relations for the intersection numbers of the
Gromov-Witten invariants twisted by the $\psi$ and $\kappa$ classes by using these
presentations and the restriction properties of these cohomology classes.
These recursion relations can be encoded in terms of a potential for these
intersection numbers satisfying a system of differential equations called the
topological recursion relations. We also prove two other equations, the
analogs of the puncture and dilaton equations, which do not use the
presentation of the $\psi$ and $\kappa$ classes and show that all of these
intersection numbers are completely determined by the Gromov-Witten
invariants of $V$.  

Consider the special case where there are no $\kappa$ classes. These
topological recursion relations were originally written down by Witten
\cite{W} in the context of topological gravity coupled to topological matter
before the moduli space of stable maps were even defined! A proof of Witten's
equations has been announced in \cite{KM2} to appear in \cite{Ma2}. Some work
toward this direction has also been done in \cite{Pa3}. (In the setting of
symplectic geometry, these equations were proven in \cite{RT}.) The proof of 
the puncture and dilaton equations for the $\psi$ classes is also forthcoming
in \cite{Ma2}. 

Our work is a generalization of these results to the case where the
$\kappa$ classes are included. The $\kappa$ classes are important here
because their restriction properties manifestly yield families of
\cft\ structures on $H^\bullet(V)$ and provide the necessary
cohomology classes (with values in the tensor algebra of
$H^\bullet(V)$) on $\M_{0,n}(V)$ which define the \cft. This is in
marked contrast to the case involving only the $\psi$ classes where the
tensor algebra-valued forms are not manifest (see \cite{KK,KMK,KMZ}
for a special case).

Furthermore, both the Gromov-Witten invariants twisted by only the $\kappa$
classes and the same twisted by only the $\psi$ classes give rise to families
of genus zero \cft\ structures on $H^\bullet(V)$. The families 
corresponding to the $\kappa$ classes and those associated to the $\psi$
classes are, in some sense, dual to each other and this correspondence
extends to higher genera \cite{KK2}.

The paper is structured as follows. In the first section, we review the
geometry of the moduli space of genus $0$ stable maps $\M_{0,n}(V)$ with $n$
marked points associated to a smooth, projective, convex variety $V$. In the
second section, we develop a graphical calculus to describe cohomology
classes on $\M_{0,n}(V)$.  In the third section, we introduce the
tautological $\psi$ and $\kappa$ classes on $\M_{0,n}(V)$. We prove
restriction properties and lifting properties of these classes and obtain a 
graphical presentation of these classes in terms of boundary strata. In the
fourth section, we show that these intersection numbers endow the cohomology
ring $H^\bullet(V)$ with families of \cft\ structures. In the fifth section,
we derive differential equations satisfied by the potentials for these
intersection numbers and prove that they are completely determined by the
Gromov-Witten invariants. 

\begin{sloppypar}
\begin{ack}
We are greatful to the Max Planck Institut f\"ur Mathematik, where this work
 was initiated, for their financial support and for providing a wonderfully
 stimulating atmosphere. We would like to thank D. Abramovich,
 P.~Belorousski, E. Getzler, Yu. Manin, R. Pandharipande, and D. Zagier for
 useful conversations. We would like to continue to thank K. Belabas for his
\TeX nical assistance and for providing the music.
\end{ack}
\end{sloppypar}


\section{The Moduli Spaces of Stable Maps}
\label{moduli}

In this section, we recall properties of the moduli space of stable maps
whose intersection numbers are the main object of study in this paper. These
spaces can be viewed as a generalization of the moduli spaces of stable curves
due to Deligne-Mumford, a canonical compactification of the moduli space of
configurations of $n$ labeled, points on Riemann surfaces of genus $g$. 

To describe the moduli space of stable maps, one needs to generalize the
notion of a stable curve.

\begin{df} 
A \emph{quasi-stable curve of genus $g$ with $n$ marked points}, denoted by
$(C ; x_1, x_2, \ldots, x_n)$, is a connected, complex projective, reduced,
curve with (at most) double points of arithmetic genus $g$ together with
distinct points $x_1, x_2, \ldots, x_n$ in $C$ away from the double points.
A \emph{special point} of $C$ is either a marked point or a double point.
\end{df}

\begin{df}
Let $V$ be a smooth, projective variety and $\beta$ in $H_2(V,\nz)$. The
\emph{moduli space of stable maps} $\M_{g,n}(V,\beta)\,:=\,\{\,[C,f;
x_1,x_2,\ldots, x_n]\,\}$ consists of isomorphism classes of the following
data: a quasi-stable curve of genus $g$ with $n$ marked points $(C; x_1, x_2,
\ldots, x_n)$ together with a morphism  $f\,:\,C\,\to\,V$ such that $\beta =
[\,f(C)\,]$ satisfying the stability condition that each irreducible
component which is mapped to a point in $V$ musts be stable in the usual
sense, {\sl i.e.}\ each irreducible component minus its special points must
have negative Euler characteristic.
\end{df}

Not all $\beta \in H_2(V,\nz)$ are represented by such curves in $V$. Let 
$B(V)$ be the semi-subgroup of $H_2(V,\nz)$ generated by the images of
the curves in $V$.


\begin{df}
The {\it evaluation morphism} associated to the $i^{\text{th}}$ marked
point $ev_i: \M_{g,n}(V,\beta) \to V$ takes the point $[C,f;
x_1,x_2,\ldots, x_n]$ to $f(x_i)$ for all $i=1,\ldots,n$. The {\it
stabilization morphism} $st: \M_{g,n}(V,\beta) \to \M_{g,n}$, defined when
$2g-2+n>0$, takes the point $[C,f; x_1,x_2,\ldots, x_n]$ to $[C;
x_1,x_2,\ldots, x_n]^{\text{st}}$. The latter means that one contracts all
components that became unstable after forgetting $f$. Let  $\pi:
\M_{g,n+1}(V,\beta) \to \M_{g,n}(V,\beta)$ be the morphism which forgets 
the $i^{\text{th}}$ point followed by stabilization, {\sl i.e.}\ $[C,f;
x_1,x_2,\ldots, x_{n+1}]\,\mapsto\,[C,f; x_1,\ldots, \widehat{x_i},\ldots,
x_{n+1}]^{\text{st}}$. This morphism is the universal curve over
$\M_{g,n}(V,\beta)$ 
\cite{Be,BMa}. 
\end{df}

\begin{thm}{\cite{BF},\cite{FP}} Let $V$ be a smooth, projective variety and
$\beta$ in $H_2(V,\nz)$. The moduli spaces $\M_{g,n}(V,\beta)$ exist and are
complete, Deligne-Mumford stacks.
\end{thm}

The moduli space of stable maps becomes the moduli space of curves by
choosing the variety $V$ to be a point. The moduli space of curves $\M_{g,n}$
is a smooth stack (orbifolds in the algebraic category) for all genera.
In general, $\M_{g,n}(V,\beta)$ is highly singular. A useful criterion by
which one can obtain nice moduli spaces, at least in genus zero, is that $V$
be convex. 

\begin{df}
A smooth, projective variety $V$ is said to be \emph{convex} if for all
morphisms $f\,:\,\CP{1}\,\to\,V$, the cohomology $H^1(\CP{1},f^*T_V)$
vanishes where $T_V$ is the (holomorphic) tangent bundle.
\end{df}

\begin{thm}{\cite{BMa},\cite{FP}}
\label{convex}
Let $V$ be a smooth, projective, convex variety of dimension $D$ then
the space $\M_{0,n}(V,\beta)$ is a locally normal, projective variety of pure
(complex) dimension $D+n-3+\int_\beta\,c_1(T_V)$ and can have, at worst,
orbifold singularities. In addition, the divisor corresponding to the
locus of the singular curves is a normal crossing divisor up to a
quotient by a finite group. 
\end{thm}

\begin{df} Let $V$ be a smooth, projective, convex variety then the {\sl
(genus zero) Gromov-Witten invariants of $V$}\, are maps
$(H^\bullet(V))^{\otimes n}\,\to\,\nc$ given by 
\[
\gamma_1\,\otimes\,\gamma_2\,\otimes\,\ldots\,\otimes\,\gamma_n\,\mapsto\, 
\left<\,\textrm{ev}_1^*\gamma_1\,\textrm{ev}_2^*\gamma_2\,\ldots\,
\textrm{ev}_n^*\gamma_n \,\right>_{\beta} \]
where 
\[ \left<\,\textrm{ev}_1^*\gamma_1\,\textrm{ev}_2^*\gamma_2\,\ldots\,
\textrm{ev}_n^*\gamma_n \,\right>_{\beta} \,:=\,\int_{\M_{0,n}(V,\beta)}\,
\textrm{ev}_1^*\gamma_1\,\textrm{ev}_2^*\gamma_2\,\ldots\,
\textrm{ev}_n^*\gamma_n
\]
where $\gamma_1,\gamma_2, \ldots, \gamma_n$ are elements in $H^\bullet(V)$,
$n$ is the number of marked points, $\beta$ belongs to $B(V)$, and the
product is understood to be the cup product.
\end{df}

Gromov-Witten invariants can also be defined in higher genera and for
certain nonconvex varieties provided that one uses the proper notion of
fundamental class (the so-called {\sl virtual fundamental class}
\cite{LT,BF,BMa}) of the moduli spaces of stable maps. 

The moduli space $\M_{0,n}(V,\beta)$ for $V$ a smooth, projective, convex
variety is endowed with a canonical stratification by strata which are
indexed by certain decorated trees called stable trees.

Let us introduce some notation regarding graphs (and trees). All graphs we
shall consider are connected and have vertices with a valence greater than or
equal to zero.  Given such a tree $\Gamma$, let $V(\Gamma)$ be its set of
vertices, $E(\Gamma)$ be its  set of edges, and $S(\Gamma)$ be its set of
tails.  Each edge has two endpoints belonging to $V(\Gamma)$ which are
allowed to be the same. Each tail has only one endpoint. For all $v\in
V(\Gamma)$,  let $n(v)$ be the number of half-edges emanating from $v$. Each
edge gives rise to two half-edges, and each tail to one half-edge.

A {\it stable graph} consists of a quadruple $(\Gamma, g, \beta,
\mu)$, where $\Gamma$ is a connected graph as above, $g:V(\Gamma) \to
\nz_{\ge 0}$, $\beta:V(\Gamma) \to B(V)$, and $\mu$ is a bijection
between $S(\Gamma)$ and a given set $I$. We also require that for
each vertex $v$ the stability condition be satisfied: if $\beta(v)=0$,
then $2g(v)-2+n(v)>0$. We define the {\it genus, $g(\Gamma)$, of $\Gamma$}
to be $b_1(\Gamma) + \sum_{v\in V(\Gamma)} g(v)$, where $b_1$ is the first
Betti number. We also define the {\it degree, $\beta(\Gamma)$, of $\Gamma$} to
be $\sum_{v\in V(\Gamma)} \beta(v)$. A {\it stable tree} is a stable graph of
genus zero.

To each stable map $[C, f; x_1, \ldots, x_n] \in \M_{0,n}(V, \beta)$
one can associate a stable tree, called its {\it dual tree}, by
collapsing each irreducible component to a point forming a vertex, 
drawing a half-edge emanating from this vertex for each special point on that
irreducible component and connecting any two half-edges if the the two
components share a double point, and then decorating each vertex $v$ with the
arithmetic genus $g(v)$ of its corresponding irreducible component and with
$\beta(v)$ which is the image of fundamental class of that irreducible
component in $B(V)$. Finally, the tails are decorated with integers
$\{\,1,\ldots,n\,\}$ corresponding to the marked points on $C$.

In this way each stable tree $\Gamma$ with $n$ tails of genus $0$ and
degree $\beta$ determines a closed sub-orbifold $\M_\Gamma$ of
$\M_{0,n}(V,\beta)$ by taking the closure of the locus of stable maps
whose dual tree is $\Gamma$ \cite{FP}. Each $\M_\Gamma$ is a smooth
stack (orbifold) of complex codimension $|E(\Gamma)|$ in $\M_{0,n}
(V,\beta)$. We denote the orbifold fundamental class of $\M_\Gamma$ by
$[\M_\Gamma]$. According to \cite{Kol}, the stratification defined
above is finite. 

The stable trees $\Gamma$ play a further role since it indicates how
each stratum $\M_\Gamma$ can be expressed as a fibered product of
fundamental classes of the moduli spaces quotiented by the action of
the automorphism group of the stable tree. More precisely, if $v\in
V(\Gamma)$ then we denote by $\M(v)$ the moduli space
$\M_{0,n(v)}(V,\beta(v))$. Each edge of $\Gamma$ joining vertices
$v_1$ and $v_2$ determines evaluation morphisms $\M(v_1) \to V$ and
$\M(v_2) \to V$. The stratum $\M_\Gamma$ is isomorphic as a stack
(orbifold) to the fibered product of $\M(v)$, $v\in V(\Gamma)$ over
$V^{|E(\Gamma)|}$ with respect to the evaluation morphisms
corresponding to each edge quotiented by the obvious action of
$\Aut(\Gamma)$, the automorphism group of $\Gamma$. We denote by
$\fprod_{v\in V(\Gamma)} \M(v)$ the fibered product defined above, by
$\rho_\Gamma$ the composition
\[
\fprod_{v\in V(\Gamma)} \M(v) \to \M_\Gamma \to \M_{0,n}(V,
\beta), 
\]
and by $j_\Gamma$ the inclusion 
\[
\fprod_{v\in V(\Gamma)} \M(v) \,\to\, \prod_{v\in V(\Gamma)}\,
\M(v).  
\]


\section{Presentation of Cohomology Classes via Graphs}
\label{gcalc}

In this section we introduce a graphical notation for certain
cohomology classes on the moduli spaces of stable maps which uses the
graphical presentation of the boundary strata. Here we generalize the
notation introduced in \cite[2.2]{KK}.

Throughout the rest of this paper we assume that $g=0$ and that $V$
is a convex, smooth, projective variety.

We start with explaining how to push forward and pull back the
cohomology classes represented by graphs under the universal curve
morphism. Let $\pi: \M_{0,n+1}(V,\beta) \to \M_{0,n}(V,\beta)$ be
defined by forgetting the $i^{\text{th}}$ marked point. First we
explain the push forward. Suppose that $\Gamma'$ represents a stratum
in $\M_{0,n+1}(V,\beta)$. If removing the tail labeled by $i$ from
$\Gamma'$ does not destabilize it, then $\pi_* ([\M_\Gamma'])=0$. On the
other hand, if removing the tail labeled $i$ does destabilize
$\Gamma'$, then
\[
\pi_* ([\M_{\Gamma'}]) = \frac{|\Aut\Gamma|}{|\Aut \Gamma'|}
[\M_\Gamma],
\] 
where $\Gamma$ is obtained by stabilization \cite{Mu}.

Next we explain how to pull back the cohomology classes represented by
graphs. Let $\M_\Gamma$ be a stratum in $\M_{0,n}(V,\beta)$. The pull
back of its fundamental class is a subvariety of $\M_{0,n+1}(V,
\beta)$ corresponding to the sum of $|V(\Gamma)|$ graphs with rational
coefficients. Each of these graphs is obtained by attaching a tail
numbered $i$ to a vertex of $\Gamma$ \cite{Mu}, call this graph
$\Gamma'$, and the corresponding coefficient is $\frac{|\Aut
\Gamma'|}{|\Aut \Gamma|}$.

These rules for pushing forward and pulling back cohomology classes on
the moduli space of stable maps are analogous to the rules for the
moduli space of curves.

In the sequel we will need trees whose vertices are decorated with
more than its degree. In addition, we will decorate the vertices with
cohomology classes. Let $\Gamma$ be a stable tree that determines a
subvariety of $\M_{0,n}(V,\beta)$, and let $\gamma_v \in H^\bullet
(\M(v))$, $v\in V(\Gamma)$. We denote the cohomology class 
\[
\frac{1}{|\Aut \Gamma|} \rho_{\Gamma *} j_\Gamma^* (\otimes_v \gamma_v)
\]
by the picture of $\Gamma$ where each vertex $v$ is in addition
labeled by the cohomology class $\gamma_v$. We omit this label if
$\gamma_v=1$, i.e., $\gamma_v$ is the orbifold fundamental class of
$\M(v)$. In particular, $[\M_\Gamma]$ is represented by $\Gamma$ with
no additional labels attached.

Before we prove the lemma below that describes the push forward of the
cohomology classes represented by decorated stable trees we would like
to recall some properties of the fibered products. Let $X$ and $Y$ be
two smooth varieties or orbifolds with morphisms $X \to V$, $Y \to V$ to a
smooth variety $V$. Then the fibered product $X \times_V Y$ is a
subspace of $X \times Y$, and we denote the corresponding inclusion by
$j$. The cohomology $H^\bullet (X \times_V Y)$ is isomorphic to
$H^\bullet (X) \otimes_{H^\bullet(V)} H^\bullet (Y)$, and the
restriction $j^*$ is surjective and the push forward $j_*$ is
injective. The composition $j_* j^* (\alpha)$ is equal to $\Delta
\alpha$, where $\alpha \in H^\bullet(X\times Y)$, and $\Delta= j_*(1)$
is the Thom class. It is the pull back of the Thom class of $V$ in
$V\times V$. 

One can easily generalize the discussion above to the fibered products
of more than two spaces. If $j$ denotes the inclusion of the fibered
product into the direct product, then $j^*$ is surjective, $j_*$ is
injective, and the composition $j_* j^*$ is the multiplication by the
pull back of the Thom class of $V^k$ in $V^n$, where the inclusion of
$V^k$ into $V^n$ is determined by the fibered product. 

\begin{lm}
\label{push}
Let $\pi: \M_{0,n+1} (V,\beta) \to \M_{0,n}(V,\beta)$ be the universal
curve determined by forgetting the tail labeled $i$. Let $\Gamma'$ be
a tree that determines a subvariety in $\M_{0,n+1} (V,\beta)$ and
$\Gamma$ be the tree obtained from $\Gamma'$ by removing the tail
labeled $i$. We assume that removing the tail labeled $i$ does not
destabilize $\Gamma'$, and for each vertex $v'$ of $\Gamma'$ we denote
by $v$ the corresponding vertex of $\Gamma$. Let $w'$ be the vertex
from which the tail numbered $i$ emanates, and let $\pi(w): \M(w')\to
\M(w)$ be the universal curve over $\M(w)$ which forgets $i$.

If $\gamma' \in H^\bullet (\M_{0,n+1} (V,\beta))$ is represented by
$\Gamma'$ decorated with the cohomology classes $\gamma'_{v'}$, $v'
\in V(\Gamma')$, then its push forward under $\pi$ is the cohomology
class $\gamma$ represented by the multiple $\frac{|\Aut \Gamma|}{|\Aut
\Gamma'|}$ of $\Gamma$ decorated with the cohomology classes
$\gamma_v$, where $\gamma_v = \gamma'_{v'}$ when $v \ne w$, and
$\gamma_w = \pi(w)_* (\gamma_{w'})$.
\end{lm}

\begin{proof}
For each vertex $v'$ of $\Gamma'$ there exists the corresponding
vertex $v$ of $\Gamma$, and $\M(v')$ is isomorphic to $\M(v)$. These
isomorphisms can be chosen so that they respect the fibered products. We
denote by $X$ the product $\prod_{v \ne w} \M(v)$ and by $X_V$ the
fibered product $\fprod_{v\ne w} \M(v)$. Consider the following
commutative diagram where the horizontal arrows are natural
morphisms, and two left vertical morphisms are induced by $\pi (w)$:
\[
\begin{CD}
X \times \M(w') @<j_{\Gamma'}<< X_V\times_V \M(w') @>\rho_{\Gamma'}>> \M_{0,n+1}(V,\beta) \\
@VV \pit V @VV \pit V @VV\pi V \\
X \times \M(w) @<j_\Gamma<<  X_V\times_V \M(w) @>\rho_\Gamma >> \M_{0,n}(V,\beta). \\
\end{CD}
\]
In order to prove the statement of the lemma one needs to show that 
\[
\pi_*\,\rho_{\Gamma' *}\,j_{\Gamma'}^{*}\,=\,\rho_{\Gamma *}\,j_\Gamma^*\,
\pit_*.
\]
This reduces to showing that $\pit_*\,j_{\Gamma'}^{*}\, = \,j_\Gamma^*
\pit_*$. The latter follows from \cite[6.2]{Fu}.
\end{proof}

\begin{rem}
If in the lemma above one assumes that $\Gamma'$ destabilizes after
forgetting the tail labeled by $i$, then $\M(w')$ is isomorphic to a
point. If the degree of $\gamma_{w'}$ is greater than $0$ attached to
$w'$, then the decorated graph represents $0$ in cohomology. If
$\gamma_{w'}=1$, then one pushes the corresponding cohomology class
down as in the case without decorations --- by stabilizing and dividing by
the orders of the appropriate automorphism groups.
\end{rem}

For example, consider the cohomology class $\gamma$ in $\M_{0,5}(V,\beta)$
(where $\beta$ is assumed to be nonzero) represented by the decorated stable
tree 
\[
\lift{1500}{$\gamma\,:=\,$} \ \begin{picture}(0,0)%
\includegraphics{./pic/gammap.pstex}%
\end{picture}%
\setlength{\unitlength}{0.00033300in}%
\begingroup\makeatletter\ifx\SetFigFont\undefined
\def\x#1#2#3#4#5#6#7\relax{\def\x{#1#2#3#4#5#6}}%
\expandafter\x\fmtname xxxxxx\relax \def\y{splain}%
\ifx\x\y   
\gdef\SetFigFont#1#2#3{%
  \ifnum #1<17\tiny\else \ifnum #1<20\small\else
  \ifnum #1<24\normalsize\else \ifnum #1<29\large\else
  \ifnum #1<34\Large\else \ifnum #1<41\LARGE\else
     \huge\fi\fi\fi\fi\fi\fi
  \csname #3\endcsname}%
\else
\gdef\SetFigFont#1#2#3{\begingroup
  \count@#1\relax \ifnum 25<\count@\count@25\fi
  \def\x{\endgroup\@setsize\SetFigFont{#2pt}}%
  \expandafter\x
    \csname \romannumeral\the\count@ pt\expandafter\endcsname
    \csname @\romannumeral\the\count@ pt\endcsname
  \csname #3\endcsname}%
\fi
\fi\endgroup
\begin{picture}(4275,1689)(1351,-2059)
\put(2551,-1636){\makebox(0,0)[lb]{\smash{\SetFigFont{8}{9.6}{rm}$\beta_1$}}}
\put(4351,-1636){\makebox(0,0)[lb]{\smash{\SetFigFont{8}{9.6}{rm}$\beta_2$}}}
\put(5626,-586){\makebox(0,0)[lb]{\smash{\SetFigFont{8}{9.6}{rm}$\phantom{1}$}}}
\put(1351,-2011){\makebox(0,0)[lb]{\smash{\SetFigFont{8}{9.6}{rm}$3$}}}
\put(5551,-811){\makebox(0,0)[lb]{\smash{\SetFigFont{8}{9.6}{rm}$4$}}}
\put(5551,-2011){\makebox(0,0)[lb]{\smash{\SetFigFont{8}{9.6}{rm}$5$}}}
\put(2026,-586){\makebox(0,0)[lb]{\smash{\SetFigFont{8}{9.6}{rm}$\phantom{1}$}}}
\put(1351,-1336){\makebox(0,0)[lb]{\smash{\SetFigFont{8}{9.6}{rm}$2$}}}
\put(1351,-736){\makebox(0,0)[lb]{\smash{\SetFigFont{8}{9.6}{rm}$1$}}}
\put(4351,-961){\makebox(0,0)[lb]{\smash{\SetFigFont{8}{9.6}{rm}$\gamma'$}}}
\put(2701,-961){\makebox(0,0)[lb]{\smash{\SetFigFont{8}{9.6}{rm}$\gamma''$}}}
\end{picture}

\]
where it is understood that $\beta_1+\beta_2\,=\,\beta$ and $\gamma'$ is a
cohomology class on the space $\M_{0,3}(V,\beta_2)$ associated to the right
vertex and $\gamma''$ is a cohomology class on the space $\M_{0,4}(V,\beta_1)$
associated to the left vertex.   Consider the morphism
$\pi\,:\,\M_{0,5}(V,\beta)\,\to\,\M_{0,4}(V,\beta)$ which forgets the
$5^{\textrm{th}}$ marked point. If $\beta_2$ is nonzero then
\[
\lift{1500}{$\pi_*(\gamma)\,:=\,$} \ \begin{picture}(0,0)%
\includegraphics{./pic/gamma_1.pstex}%
\end{picture}%
\setlength{\unitlength}{0.00033300in}%
\begingroup\makeatletter\ifx\SetFigFont\undefined
\def\x#1#2#3#4#5#6#7\relax{\def\x{#1#2#3#4#5#6}}%
\expandafter\x\fmtname xxxxxx\relax \def\y{splain}%
\ifx\x\y   
\gdef\SetFigFont#1#2#3{%
  \ifnum #1<17\tiny\else \ifnum #1<20\small\else
  \ifnum #1<24\normalsize\else \ifnum #1<29\large\else
  \ifnum #1<34\Large\else \ifnum #1<41\LARGE\else
     \huge\fi\fi\fi\fi\fi\fi
  \csname #3\endcsname}%
\else
\gdef\SetFigFont#1#2#3{\begingroup
  \count@#1\relax \ifnum 25<\count@\count@25\fi
  \def\x{\endgroup\@setsize\SetFigFont{#2pt}}%
  \expandafter\x
    \csname \romannumeral\the\count@ pt\expandafter\endcsname
    \csname @\romannumeral\the\count@ pt\endcsname
  \csname #3\endcsname}%
\fi
\fi\endgroup
\begin{picture}(4275,1689)(1351,-2059)
\put(2551,-1636){\makebox(0,0)[lb]{\smash{\SetFigFont{8}{9.6}{rm}$\beta_1$}}}
\put(4351,-1636){\makebox(0,0)[lb]{\smash{\SetFigFont{8}{9.6}{rm}$\beta_2$}}}
\put(5626,-586){\makebox(0,0)[lb]{\smash{\SetFigFont{8}{9.6}{rm}$\phantom{1}$}}}
\put(1351,-2011){\makebox(0,0)[lb]{\smash{\SetFigFont{8}{9.6}{rm}$3$}}}
\put(5551,-811){\makebox(0,0)[lb]{\smash{\SetFigFont{8}{9.6}{rm}$4$}}}
\put(2026,-586){\makebox(0,0)[lb]{\smash{\SetFigFont{8}{9.6}{rm}$\phantom{1}$}}}
\put(1351,-1336){\makebox(0,0)[lb]{\smash{\SetFigFont{8}{9.6}{rm}$2$}}}
\put(1351,-736){\makebox(0,0)[lb]{\smash{\SetFigFont{8}{9.6}{rm}$1$}}}
\put(3826,-961){\makebox(0,0)[lb]{\smash{\SetFigFont{8}{9.6}{rm}$\tilde{\pi}_*(\gamma')$}}}
\put(2626,-961){\makebox(0,0)[lb]{\smash{\SetFigFont{8}{9.6}{rm}$\gamma''$}}}
\end{picture}

\]
where $\pit\,:\,\M_{0,3}(V,\beta_2)\,\to\, \M_{0,2}(V,\beta_2)$ is the
projection forgetting the point labeled by $5$ associated to the right
vertex. On the other hand, if $\beta_2\,=\,0$ then the right hand side of the
previous equation is now unstable. If $\gamma'$ has nonzero dimension
then $\pi_*(\gamma)$ vanishes. If $\gamma'$ is the unit element then
\[
\lift{1500}{$\pi_*(\gamma)\,:=\,$} \ \begin{picture}(0,0)%
\includegraphics{./pic/gamma_2.pstex}%
\end{picture}%
\setlength{\unitlength}{0.00033300in}%
\begingroup\makeatletter\ifx\SetFigFont\undefined
\def\x#1#2#3#4#5#6#7\relax{\def\x{#1#2#3#4#5#6}}%
\expandafter\x\fmtname xxxxxx\relax \def\y{splain}%
\ifx\x\y   
\gdef\SetFigFont#1#2#3{%
  \ifnum #1<17\tiny\else \ifnum #1<20\small\else
  \ifnum #1<24\normalsize\else \ifnum #1<29\large\else
  \ifnum #1<34\Large\else \ifnum #1<41\LARGE\else
     \huge\fi\fi\fi\fi\fi\fi
  \csname #3\endcsname}%
\else
\gdef\SetFigFont#1#2#3{\begingroup
  \count@#1\relax \ifnum 25<\count@\count@25\fi
  \def\x{\endgroup\@setsize\SetFigFont{#2pt}}%
  \expandafter\x
    \csname \romannumeral\the\count@ pt\expandafter\endcsname
    \csname @\romannumeral\the\count@ pt\endcsname
  \csname #3\endcsname}%
\fi
\fi\endgroup
\begin{picture}(4275,1689)(1351,-2059)
\put(5626,-586){\makebox(0,0)[lb]{\smash{\SetFigFont{8}{9.6}{rm}$\phantom{1}$}}}
\put(1351,-2011){\makebox(0,0)[lb]{\smash{\SetFigFont{8}{9.6}{rm}$3$}}}
\put(2026,-586){\makebox(0,0)[lb]{\smash{\SetFigFont{8}{9.6}{rm}$\phantom{1}$}}}
\put(1351,-1336){\makebox(0,0)[lb]{\smash{\SetFigFont{8}{9.6}{rm}$2$}}}
\put(1351,-736){\makebox(0,0)[lb]{\smash{\SetFigFont{8}{9.6}{rm}$1$}}}
\put(2551,-1636){\makebox(0,0)[lb]{\smash{\SetFigFont{8}{9.6}{rm}$\beta$}}}
\put(3826,-1336){\makebox(0,0)[lb]{\smash{\SetFigFont{8}{9.6}{rm}$4$}}}
\put(2626,-961){\makebox(0,0)[lb]{\smash{\SetFigFont{8}{9.6}{rm}$\gamma''$}}}
\end{picture}

\]
In all of these cases, the automorphism group is trivial so the prefactor
does not arise.


\section{Tautological Classes}
\label{taut}

In this section we introduce the tautological $\psi$ and $\kappa$
classes on the moduli spaces of stable maps which generalize the
corresponding tautological classes on the moduli spaces of stable
curves. We will also discuss the properties of the pull back of the
$\psi$ classes under the stabilization morphism, and derive the
restriction properties of the tautological classes to the boundary
strata. 

Let $\pi: \M_{0,n+1}(V, \beta) \to \M_{0,n} (V,\beta)$ be the
universal curve. We assume that $\pi$ ``forgets'' the
$(n+1)^{\text{st}}$ marked point. The morphism $\pi$ has $n$ canonical
sections $\sigma_1, \ldots, \sigma_n$. We denote their images by
$D_1, \ldots, D_n$. If $g=0$ and $V$ is convex each $D_i$ is a
divisor. We denote by $\omega$ the relative dualizing sheaf of $\pi$.
It is an invertible sheaf, that is, it determines a line bundle on
$\M_{0,n+1}(V, \beta)$ \cite{BMa}. 

\begin{df}
For each $i=1,\ldots, n$ the {\it tautological line bundle} $\LL_i$ on
$\M_{0,n}(V,\beta)$ is $\sigma_i^* \omega$. The {\it tautological
class} $\psi_i \in H^2(\M_{0,n}(V,\beta))$ is equal to the first Chern
class $c_1(\LL_i)$.
\end{df}

One can also pull back cohomology classes from $V$ to
$\M_{0,n}(V,\beta)$ using the evaluation maps. The definition of the
$\kappa$ classes involves both, powers of the $\psi$ classes and these
pull backs. 

\begin{df}
The {\it tautological class} $\kappa_{a}$ in
$H^{\bullet}(\M_{0,n}(V,\beta))\, \otimes\,H^\bullet(V)^*$ for $a\ge
-1$ is defined as follows. For each $\gamma\,\in\,H^\bullet(V)$, the
cohomology class $\kappa_{a}(\gamma)$ is the push forward of
$\psi_{n+1}^{a+1} ev_{n+1}^* (\gamma)$ with respect to the projection
$\pi: \M_{0,n+1}(V,\beta) \to \M_{0,n}(V,\beta)$ which forgets the
$(n+1)^{\text{st}}$ marked point. In particular, if $\gamma$ has
definite degree $|\gamma|$ then $\kappa_{a} (\gamma)$ has degree $2a +
|\gamma|$. If $\{\, e_\alpha\, \}_{\alpha \in A}$, is a homogeneous
basis for $H^\bullet(V)$ such that $e_0$ is the identity element then
$\kappa_{a,\alpha}$ denotes the cohomology class $\kappa_a
(e_\alpha)$.
\end{df}

\begin{rem}
The class $\kappa_{-1,0}$ vanishes due to dimensional reasons. In
addition, all classes $\kappa_{-1}(\gamma)$ vanish on $\M_{0,n}(V,0)$.
\end{rem}

Our definition corresponds to the ``modified'' $\kappa$ classes
defined by Arbarello and Cornalba \cite{AC} rather than the
``classical'' $\kappa$ classes defined by Mumford \cite{Mu}. It is
easy to see that $\psi_{n+1} = c_1 (\omega (D_1 + \ldots + D_n))$, and
therefore one can also define the $\kappa$ classes using the relative
dualizing sheaf.

Next we will show that the $\psi$ classes restrict to the strata in the
expected manner.  For a more general treatment, we again refer to
\cite{Ma2} (see also \cite{Pa}).

\begin{lm}
\label{psi:restr}
Let $\Gamma$ be a stable tree which determines a stratum in
$\M_{0,n}(V,\beta)$. Suppose that the tail labeled $i$ is attached
to $w\in V(\Gamma)$. Denote the class $\psi_i$ on $\M_{0,n}(V,\beta)$
by $\psi_i$, and on $\M(w)$ by $\psi'_i$. Then $\rho_\Gamma^* (\psi_i) =
j_\Gamma^* \psi'_i$.
\end{lm}

\begin{proof}
Let $\pi: \M_{0,n+1} (V,\beta) \to \M_{0,n}(V,\beta)$ and $\pi (w):
\M(w') \to \M(w)$ be the universal curves forgetting the marked point
labeled $n+1$, and $\sigma_i$ be the corresponding sections which we
denote by the same symbol. We use the same notation as in the proof of
Lem.~\ref{push}. Consider the commutative diagram
\[
\begin{CD}
X \times \M(w') @<j_{\Gamma'}<< X_V\times_V \M(w') @>\rho_{\Gamma'}>> \M_{0,n+1}(V,\beta) \\
@VV \pit V @VV\pit V @VV\pi V \\
X \times \M(w) @<j_\Gamma<<  X_V\times_V \M(w) @>\rho_\Gamma >> \M_{0,n}(V,\beta). \\
\end{CD}
\]

Let us denote by $\omega$ with a subscript the relative dualizing
sheaf of the morphism corresponding to the subscript. It is clear that
$\rho_1^* \omega_\pi = \omega_{1\times_V \pit} = j_1^* \omega_{1\times
\pit}$ when restricted to the image of $\sigma_i$. It follows that
\[
\rho_\Gamma^* \psi_i = \sigma_i^* \rho_{\Gamma'}^* (\omega_\pi) = 
\sigma_i^* j_{\Gamma'}^* (\omega_{\pit}) = j_\Gamma^* \psi_i^{'}. 
\]
\end{proof}

It is easy to verify that the pull backs of the cohomology classes on
$V$ under the evaluation maps restrict in the same manner.

\begin{lm}
\label{kappa:restr}
Let $\Gamma$ be a stable tree which determines a stratum in
$\M_{0,n}(V,\beta)$. Denote the class $\kappa_{a,\alpha}$ on $\M_{0,n}
(V,\beta)$ (resp. $\M(v)$, where $v\in V(\Gamma))$ by $\kappa$ (resp.
$\kappa (v))$. Then 
\[
\rho_\Gamma^* (\kappa) = j_\Gamma^* \sum_{v\in V(\Gamma)} \kappa(v).
\]
\end{lm}

\begin{proof}
Consider the set $\mathcal{G}$ of cardinality $|V(\Gamma)|$ whose
elements are graphs $\Gamma'$ each of which is obtained from $\Gamma$ by
attaching a tail labeled $n+1$ to one of the vertices of $\Gamma$.
Each $\Gamma' \in \mathcal{G}$ determines a stratum of $\M_{0,n+1}
(V,\beta)$, and there are natural morphisms 
\[
\prod_{v\in V(\Gamma')} \M(v) \to \prod_{v\in V(\Gamma)} \M(v),
\]
and the induced one on the fibered products. We denote both of them by
$\pit$. 

Consider the following commutative diagram

\[
\begin{CD}
\coprod_{\Gamma'\in \mathcal{G}} \prod_{v\in \Gamma'} \M(v) @< j_{\Gamma'}<<
\coprod_{\Gamma' \in \mathcal{G}} \fprod_{v\in \Gamma'} \M(v) @> \rho_{\Gamma'}>>
\M_{0,n+1} (V,\beta) \\
@VV \pit V @VV \pit V @VV\pi V \\
\prod_{v\in \Gamma} \M(v) @< j_{\Gamma}<<
\fprod_{v\in \Gamma} \M(v) @> \rho_{\Gamma}>>
\M_{0,n} (V,\beta). \\
\end{CD}
\]

Since $j_{\Gamma *}$ is injective it suffices to show that $j_{\Gamma
*} \rho_\Gamma^* \kappa = \Delta \sum_v \kappa(v)$, where $\Delta$ is the
Thom class of $\fprod_v \M(v)$ in $\prod_v \M(v)$.

It suffices to show that 
\begin{equation}
\label{push:pull}
\pit_* j_{\Gamma' *} \rho_{\Gamma'}^* = 
j_{\Gamma *} \rho_\Gamma^* \pi_*.
\end{equation}
Indeed, 
\begin{multline*}
j_{\Gamma *} \rho_\Gamma^* \kappa = 
j_{\Gamma *} \rho_\Gamma^* \pi_* (\psi_{n+1}^{a+1} ev_{n+1}(e_\alpha))  =
\pit_* j_{\Gamma' *} \rho_{\Gamma'}^* (\psi_{n+1}^{a+1}
ev_{n+1}(e_\alpha)) \\ =
\pit_* j_{\Gamma' *} j_{\Gamma'}^* (\psi_{n+1}^{a+1}
ev_{n+1}(e_\alpha)) = 
\pit_* (\Deltat \psi_{n+1}^{a+1}
ev_{n+1}(e_\alpha)) =
\Delta \sum_v \kappa(v). 
\end{multline*}
Here we have used \eqref{push:pull} for the second equality,
Lem.~\ref{psi:restr} for the third equality. We denoted by $\Deltat$
the Thom class of $\coprod_{\Gamma' \in \mathcal{G}} \fprod_{v\in
\Gamma'} \M(v)$ in $\coprod_{\Gamma'\in \mathcal{G}} \prod_{v\in
\Gamma'} \M(v)$, and have used that it is the pull back of $\Delta$
under the map $\pit$.

In order to show \eqref{push:pull} it suffices to prove that $\pit_*
\rho_{\Gamma'}^* = \rho_\Gamma^* \pi_*$. This follows from the fact
that the right square is the normalization of a pull back square and
\cite[6.2]{Fu}.
\end{proof}

Our next goal is to express the $\psi$ classes via the boundary
strata. Since the action of the symmetric group on the moduli spaces of
stable maps interchanges the $\psi$ classes it is enough to give a
presentation of $\psi_1$ via the boundary strata. We will denote
$\psi_1$ on $\M_{0,n}(V,\beta)$ by $\psi_{(n,\beta)}$. We also assume
that if a graph is not stable, then it represents the zero cohomology
class. 

\begin{lm}
\label{psi:pres}
If $n\ge 3$, $a\ge 1$, then the following holds in $H^\bullet
(\M_{0,n} (V,\beta))$
\[
\lift{1500}{$\psi_{(n,\beta)}^a = 
\displaystyle{ 
\sum_{\substack{\beta_1+\beta_2 = \beta \\ I \sqcup J =
[n-2] \\ 1 \in J}}}$}
\ \begin{picture}(0,0)%
\includegraphics{./pic/psi.pstex}%
\end{picture}%
\setlength{\unitlength}{0.00033300in}%
\begingroup\makeatletter\ifx\SetFigFont\undefined
\def\x#1#2#3#4#5#6#7\relax{\def\x{#1#2#3#4#5#6}}%
\expandafter\x\fmtname xxxxxx\relax \def\y{splain}%
\ifx\x\y   
\gdef\SetFigFont#1#2#3{%
  \ifnum #1<17\tiny\else \ifnum #1<20\small\else
  \ifnum #1<24\normalsize\else \ifnum #1<29\large\else
  \ifnum #1<34\Large\else \ifnum #1<41\LARGE\else
     \huge\fi\fi\fi\fi\fi\fi
  \csname #3\endcsname}%
\else
\gdef\SetFigFont#1#2#3{\begingroup
  \count@#1\relax \ifnum 25<\count@\count@25\fi
  \def\x{\endgroup\@setsize\SetFigFont{#2pt}}%
  \expandafter\x
    \csname \romannumeral\the\count@ pt\expandafter\endcsname
    \csname @\romannumeral\the\count@ pt\endcsname
  \csname #3\endcsname}%
\fi
\fi\endgroup
\begin{picture}(4577,1479)(1276,-1873)
\put(2026,-586){\makebox(0,0)[lb]{\smash{\SetFigFont{8}{9.6}{rm}$n\!-\!1$}}}
\put(3001,-586){\makebox(0,0)[lb]{\smash{\SetFigFont{8}{9.6}{rm}$n$}}}
\put(1276,-1261){\makebox(0,0)[lb]{\smash{\SetFigFont{8}{9.6}{rm}$I$}}}
\put(2551,-1636){\makebox(0,0)[lb]{\smash{\SetFigFont{8}{9.6}{rm}$\beta_1$}}}
\put(5853,-1263){\makebox(0,0)[lb]{\smash{\SetFigFont{8}{9.6}{rm}$J$}}}
\put(5626,-586){\makebox(0,0)[lb]{\smash{\SetFigFont{8}{9.6}{rm}$1$}}}
\put(4351,-1636){\makebox(0,0)[lb]{\smash{\SetFigFont{8}{9.6}{rm}$\beta_2$}}}
\put(4351,-886){\makebox(0,0)[lb]{\smash{\SetFigFont{8}{9.6}{rm}$\psi^{a-1}$}}}
\end{picture}

\]
\end{lm}

\begin{proof}
Getzler \cite{Ge} proved the lemma for all projective manifolds in
case $a=1$ using techniques involving algebraic stacks. It is stated
in \cite{KM2} that a proof of the lemma for all projective manifolds
in case $a=1$ will appear in \cite{Ma2}. Pandharipande \cite{Pa}
informed us that he also has a proof of this lemma (unpublished). We
present a proof in case of convex $V$ that does not use stacks.

It is enough to prove the lemma for $a=1$. This and
Lem.~\ref{psi:restr} will imply the lemma for all $a\ge 1$.

First we prove that 
\begin{equation*}
\psi_{(3,\beta)} = \sum_{\beta_1 + \beta_2 = \beta}
[D_{\beta_1,\beta_2}],
\end{equation*}
where 
\begin{equation*}
\lift{1500}{$[D_{\beta_1,\beta_2}] = 
\displaystyle{ 
\sum_{\substack{\beta_1+\beta_2 = \beta \\ I \sqcup J =
[n-2] \\ 1 \in J}}}$}
\ \begin{picture}(0,0)%
\includegraphics{./pic/psi_03.pstex}%
\end{picture}%
\setlength{\unitlength}{0.00033300in}%
\begingroup\makeatletter\ifx\SetFigFont\undefined
\def\x#1#2#3#4#5#6#7\relax{\def\x{#1#2#3#4#5#6}}%
\expandafter\x\fmtname xxxxxx\relax \def\y{splain}%
\ifx\x\y   
\gdef\SetFigFont#1#2#3{%
  \ifnum #1<17\tiny\else \ifnum #1<20\small\else
  \ifnum #1<24\normalsize\else \ifnum #1<29\large\else
  \ifnum #1<34\Large\else \ifnum #1<41\LARGE\else
     \huge\fi\fi\fi\fi\fi\fi
  \csname #3\endcsname}%
\else
\gdef\SetFigFont#1#2#3{\begingroup
  \count@#1\relax \ifnum 25<\count@\count@25\fi
  \def\x{\endgroup\@setsize\SetFigFont{#2pt}}%
  \expandafter\x
    \csname \romannumeral\the\count@ pt\expandafter\endcsname
    \csname @\romannumeral\the\count@ pt\endcsname
  \csname #3\endcsname}%
\fi
\fi\endgroup
\begin{picture}(4200,1689)(1426,-2059)
\put(2551,-1636){\makebox(0,0)[lb]{\smash{\SetFigFont{8}{9.6}{rm}$\beta_1$}}}
\put(5626,-586){\makebox(0,0)[lb]{\smash{\SetFigFont{8}{9.6}{rm}$\phantom{1}$}}}
\put(4351,-1636){\makebox(0,0)[lb]{\smash{\SetFigFont{8}{9.6}{rm}$\beta_2$}}}
\put(5626,-1336){\makebox(0,0)[lb]{\smash{\SetFigFont{8}{9.6}{rm}$1$}}}
\put(2026,-586){\makebox(0,0)[lb]{\smash{\SetFigFont{8}{9.6}{rm}$\phantom{n\!-\!1}$}}}
\put(1426,-736){\makebox(0,0)[lb]{\smash{\SetFigFont{8}{9.6}{rm}$2$}}}
\put(1426,-2011){\makebox(0,0)[lb]{\smash{\SetFigFont{8}{9.6}{rm}$3$}}}
\end{picture}

\end{equation*}
and $D_{\beta_1, \beta_2}$ is the corresponding divisor. The stability
condition implies that $\beta_2 \ne 0$.

First note that the restriction of $\psi_{(3,\beta)}$ to $U:=
\M_{0,3}(V,\beta) - \bigcup D_{\beta_1,\beta_2}$ is zero. Indeed,
consider the universal map $\M_{0,4} (V,\beta) \to \M_{0,3} (V,\beta)$
restricted to $U$. Let $x\in U$, and $C_x$ be the inverse image of $x$
under the universal map. The image of the section of the universal map
determined by the first marked point lies on the irreducible component
of $C_x$ which remains stable after forgetting the map $f_x: C_x \to
V$. This shows that $\psi_{(3,\beta)}$ restricted to $U$ is the
pull-back from $\M_{0,3}$, and therefore trivial. 

It follows from the previous paragraph that 
\begin{equation} \label{psi:psi}
\psi_{(3,\beta)} = \sum_{\beta_1 + \beta_2 = \beta, \, i} 
\, k_{\beta_1, \beta_2, i}\, [D_{\beta_1, \beta_2, i}],
\end{equation}
where $D_{\beta_1, \beta_2, i}$ are (disjoint) irreducible components
of $D_{\beta_1, \beta_2}$. In order to determine the coefficients it
suffices to intersect \eqref{psi:psi} with an arbitrary complete
non-singular curve $S$ in $\M_{0,3}(V,\beta)$. Let $T \to S$ be a
family of $3$-pointed stable maps over $S$, and $\nu$ be the
corresponding morphism $S \to \M_{0,3} (V,\beta)$. We can assume that
the image of $\nu$ intersects each codimension one boundary stratum
transversally, and it is disjoint from the intersection of the
boundary strata. Then $[S] \cdot \sum [D_{\beta_1, \beta_2}] =
\sum_{\beta_1, \beta_2, i} l_{\beta_1, \beta_2, i}$, where
$l_{\beta_1, \beta_2, i}$ is the number of points in the intersection
of $S$ with $D_{\beta_1, \beta_2, i}$. The stabilization morphism,
i.e., the morphism determined by forgetting the map to $V$, determines
a morphism $T \to S \times \CP{1}$. It is the blow up at the points
corresponding to contracted rational components. The image of the
first section contains exactly those points which lie in the fibers
over $D_{\beta_1,\beta_2,i} \cap S$. It follows that $\psi_{(3,\beta)}
\cdot [S] = \sum [D_{\beta_1, \beta_2}] \cdot [S]$, and therefore
$\psi_{(3,\beta)} = \sum [D_{\beta_1, \beta_2}]$.

In order to obtain the statement of the lemma for $n\ge 4$ one can
first show that if $\pi: \M_{0,n+1}(V,\beta) \to \M_{0,n}(V,\beta)$ is
the universal map, then $\psi_{(n+1,\beta)} = \pi^* \psi_{(n,\beta)} +
[D_{1,n+1}]$, where $D_{1,n+1}$ is the image of the first section of the
universal map. This can be done by using an argument similar to the above.
(See also \cite{Kn}.) Then one need to pull back cohomology
classes represented by the boundary strata using the description from
Sec.~\ref{gcalc}. 
\end{proof}

As a corollary using Lem.~\ref{push} we express the $\kappa$ classes
via the boundary strata and the $\kappa_{-1}$ classes. We denote the
class $\kappa_{a,\alpha}$ on $\M_{0,n}(V,\beta)$ by $\kappa_{a,\alpha,
(n,\beta)}$.

\begin{crl}
\label{kappa:pres}
If $n\ge 2$, $a\ge 1$, then the following holds in $H^\bullet
(\M_{0,n}(V,\beta))$ 

\begin{equation*}
\lift{1500}{$\kappa_{a,\alpha, (n,\beta)} =
\displaystyle{ 
\sum_{\substack{\beta_1+\beta_2 = \beta \\ I \sqcup J =
[n-2]}}}$}
\ \begin{picture}(0,0)%
\includegraphics{./pic/kappa.pstex}%
\end{picture}%
\setlength{\unitlength}{0.00033300in}%
\begingroup\makeatletter\ifx\SetFigFont\undefined
\def\x#1#2#3#4#5#6#7\relax{\def\x{#1#2#3#4#5#6}}%
\expandafter\x\fmtname xxxxxx\relax \def\y{splain}%
\ifx\x\y   
\gdef\SetFigFont#1#2#3{%
  \ifnum #1<17\tiny\else \ifnum #1<20\small\else
  \ifnum #1<24\normalsize\else \ifnum #1<29\large\else
  \ifnum #1<34\Large\else \ifnum #1<41\LARGE\else
     \huge\fi\fi\fi\fi\fi\fi
  \csname #3\endcsname}%
\else
\gdef\SetFigFont#1#2#3{\begingroup
  \count@#1\relax \ifnum 25<\count@\count@25\fi
  \def\x{\endgroup\@setsize\SetFigFont{#2pt}}%
  \expandafter\x
    \csname \romannumeral\the\count@ pt\expandafter\endcsname
    \csname @\romannumeral\the\count@ pt\endcsname
  \csname #3\endcsname}%
\fi
\fi\endgroup
\begin{picture}(4577,1503)(1276,-1873)
\put(2026,-586){\makebox(0,0)[lb]{\smash{\SetFigFont{8}{9.6}{rm}$n\!-\!1$}}}
\put(3001,-586){\makebox(0,0)[lb]{\smash{\SetFigFont{8}{9.6}{rm}$n$}}}
\put(1276,-1261){\makebox(0,0)[lb]{\smash{\SetFigFont{8}{9.6}{rm}$I$}}}
\put(2551,-1636){\makebox(0,0)[lb]{\smash{\SetFigFont{8}{9.6}{rm}$\beta_1$}}}
\put(5853,-1263){\makebox(0,0)[lb]{\smash{\SetFigFont{8}{9.6}{rm}$J$}}}
\put(4351,-1636){\makebox(0,0)[lb]{\smash{\SetFigFont{8}{9.6}{rm}$\beta_2$}}}
\put(5626,-586){\makebox(0,0)[lb]{\smash{\SetFigFont{8}{9.6}{rm}$\phantom{1}$}}}
\put(3751,-886){\makebox(0,0)[lb]{\smash{\SetFigFont{8}{9.6}{rm}$\kappa_{a-1,\alpha}$}}}
\end{picture}

\end{equation*}

If $n\ge 2$, then the following holds in $H^\bullet
(\M_{0,n}(V,\beta))$ 

\begin{multline*}
\lift{1500}{$\kappa_{0,\alpha, (n,\beta)} =
\displaystyle{ 
\sum_{\substack{\beta_1+\beta_2 = \beta \\ I \sqcup J =
[n-2]}}}$}
\ \begin{picture}(0,0)%
\includegraphics{./pic/kappa_0_1.pstex}%
\end{picture}%
\setlength{\unitlength}{0.00033300in}%
\begingroup\makeatletter\ifx\SetFigFont\undefined
\def\x#1#2#3#4#5#6#7\relax{\def\x{#1#2#3#4#5#6}}%
\expandafter\x\fmtname xxxxxx\relax \def\y{splain}%
\ifx\x\y   
\gdef\SetFigFont#1#2#3{%
  \ifnum #1<17\tiny\else \ifnum #1<20\small\else
  \ifnum #1<24\normalsize\else \ifnum #1<29\large\else
  \ifnum #1<34\Large\else \ifnum #1<41\LARGE\else
     \huge\fi\fi\fi\fi\fi\fi
  \csname #3\endcsname}%
\else
\gdef\SetFigFont#1#2#3{\begingroup
  \count@#1\relax \ifnum 25<\count@\count@25\fi
  \def\x{\endgroup\@setsize\SetFigFont{#2pt}}%
  \expandafter\x
    \csname \romannumeral\the\count@ pt\expandafter\endcsname
    \csname @\romannumeral\the\count@ pt\endcsname
  \csname #3\endcsname}%
\fi
\fi\endgroup
\begin{picture}(4577,1503)(1276,-1873)
\put(2026,-586){\makebox(0,0)[lb]{\smash{\SetFigFont{8}{9.6}{rm}$n\!-\!1$}}}
\put(3001,-586){\makebox(0,0)[lb]{\smash{\SetFigFont{8}{9.6}{rm}$n$}}}
\put(1276,-1261){\makebox(0,0)[lb]{\smash{\SetFigFont{8}{9.6}{rm}$I$}}}
\put(2551,-1636){\makebox(0,0)[lb]{\smash{\SetFigFont{8}{9.6}{rm}$\beta_1$}}}
\put(5853,-1263){\makebox(0,0)[lb]{\smash{\SetFigFont{8}{9.6}{rm}$J$}}}
\put(4351,-1636){\makebox(0,0)[lb]{\smash{\SetFigFont{8}{9.6}{rm}$\beta_2$}}}
\put(5626,-586){\makebox(0,0)[lb]{\smash{\SetFigFont{8}{9.6}{rm}$\phantom{1}$}}}
\put(3751,-886){\makebox(0,0)[lb]{\smash{\SetFigFont{8}{9.6}{rm}$\kappa_{-1,\alpha}$}}}
\end{picture}
 \\
\lift{1500}{$\,+\,
\displaystyle{ \sum_{i\in [n-2]}}$} \ \begin{picture}(0,0)%
\includegraphics{./pic/kappa_0_2.pstex}%
\end{picture}%
\setlength{\unitlength}{0.00033300in}%
\begingroup\makeatletter\ifx\SetFigFont\undefined
\def\x#1#2#3#4#5#6#7\relax{\def\x{#1#2#3#4#5#6}}%
\expandafter\x\fmtname xxxxxx\relax \def\y{splain}%
\ifx\x\y   
\gdef\SetFigFont#1#2#3{%
  \ifnum #1<17\tiny\else \ifnum #1<20\small\else
  \ifnum #1<24\normalsize\else \ifnum #1<29\large\else
  \ifnum #1<34\Large\else \ifnum #1<41\LARGE\else
     \huge\fi\fi\fi\fi\fi\fi
  \csname #3\endcsname}%
\else
\gdef\SetFigFont#1#2#3{\begingroup
  \count@#1\relax \ifnum 25<\count@\count@25\fi
  \def\x{\endgroup\@setsize\SetFigFont{#2pt}}%
  \expandafter\x
    \csname \romannumeral\the\count@ pt\expandafter\endcsname
    \csname @\romannumeral\the\count@ pt\endcsname
  \csname #3\endcsname}%
\fi
\fi\endgroup
\begin{picture}(5925,1503)(1,-1873)
\put(  1,-1336){\makebox(0,0)[lb]{\smash{\SetFigFont{8}{9.6}{rm}$[n-2]-i$}}}
\put(2326,-586){\makebox(0,0)[lb]{\smash{\SetFigFont{8}{9.6}{rm}$n\!-\!1$}}}
\put(3301,-586){\makebox(0,0)[lb]{\smash{\SetFigFont{8}{9.6}{rm}$n$}}}
\put(5926,-586){\makebox(0,0)[lb]{\smash{\SetFigFont{8}{9.6}{rm}$\phantom{1}$}}}
\put(2851,-1636){\makebox(0,0)[lb]{\smash{\SetFigFont{8}{9.6}{rm}$\beta$}}}
\put(4501,-1336){\makebox(0,0)[lb]{\smash{\SetFigFont{8}{9.6}{rm}$i$}}}
\put(3226,-1036){\makebox(0,0)[lb]{\smash{\SetFigFont{8}{9.6}{rm}$\mathrm{ev}^*_i\,e_\alpha$}}}
\end{picture}

\end{multline*}
\end{crl}

It follows from the above presentation of the $\kappa$ classes that
$\kappa_{0,0} = n-2 \in H^0 (\M_{0,n} (V,\beta))$. This can be also seen
directly by integrating over the fibers of the universal curve. One
can also check that if $\gamma \in H^2(V)$, then $\kappa_{-1}(\gamma)= 
\int_\beta\,\gamma$ on $\M_{0,n} (V,\beta)$. We will need this in
Sec.~\ref{trr}.          


\section{Cohomological Field Theories}
\label{cft}

In this section, we define the notion of a cohomological field theory
and prove that the Gromov-Witten invariants twisted by the $\kappa$ classes
endows $H^\bullet(V)$ with the a family of genus zero \cft\ structures. This
is equivalent to endowing $H^\bullet(V)$ with a family of formal Frobenius
manifold structures (with nonflat identity element, in general) arising from
the Poincar\'e pairing and deformations of the cup product on $V$. Our
deformations contain quantum cohomology as a special case.

The moduli spaces of curves $\M\,:=\,\{\,\M_{g,n}\,\}$ forms a modular
operad, a higher genus generalization of an operad, \cite{GK,GK2} in the
category of smooth stacks.  The action of the symmetric group relabels the
marked points while the inclusion of strata gives rise to the composition
maps. This endows the homology groups
$H_\bullet(\M)\,:=\,\{\,H_\bullet(\M_{g,n})\,\}$ with the structure of a 
modular operad in the category of graded vector spaces.  Let $(\HH,\eta)$ be a
($\nz_2$-graded) vector space $\HH$ with an even, symmetric, nondegenerate,
bilinear form $\eta$, the endomorphism operad is given by
$\End(\HH)\,:=\{\,\End_{g,n}(\HH)\,\}$ where $\End_{g,n}(\HH)\,:=\, T^n \HH^*
$ for stable pairs $(g,n)$. 

\begin{df}
A \emph{(complete) cohomological field theory} $(\HH,\eta)$ (or \cft) is a
morphism of modular operads $H_\bullet(\M)\,\to\,\End(\HH)$, {\sl i.e.}\
$(\HH,\eta)$ is an $H_\bullet(\M)$-algebra.  A \emph{cohomological field
theory of genus $g$} is a sequence of linear maps
$H_\bullet(\M_{g',n})\,\to\,\End_{g',n}(\HH)$ for all $g'\leq g$ 
satisfying the subset of axioms of a cohomological field theory which
includes only objects of genus $g' \leq g$.
\end{df}

There is a dual description of a \cft\ in terms of cohomology classes. A
\cft\ is a pair $(\HH,\eta)$ together with a collection 
$\Omega\,:=\,\{\,\Omega_{g,n}\,\}$ where $\Omega_{g,n}$ is an even element in
$\R_{g,n}\,:=\,H^\bullet(\M_{g,n}) \otimes T^n\HH^*$ defined for stable pairs
$(g,n)$ satisfying the following (where the summation convention has been
used): 
\begin{description}
\item[i] $\Omega_{g,n}$ is invariant under the action of the symmetric
  group $S_n$.
\item[ii] For each partition of $[n]\,=\,S_1\sqcup S_2$ such that $|S_1|=n_1$
 and  $|S_2|=n_2$ and nonnegative $g_1$, $g_2$ such that $g\,=\,g_1\,+\,g_2$
 and $2g_i-2+n_i\,>\,0$ for all $i$,  consider the inclusion map 
$\rho\,:\,\M_{g_1,S_1 \sqcup *}\times\M_{g_2,S_2\sqcup
 *}\,\to\,\M_{g_1+g_2,n}$ where $*$ denotes the two marked points that are
 attached under the inclusion map. The forms satisfy the restriction property
\begin{align*}
&\rho^*\,\Omega_{g,n}(\gamma_1,\gamma_2\,\ldots,\gamma_n)\,=\, \\
&\pm\,\Omega_{g_1,n_1}((\bigotimes_{\alpha\in S_1}\,\gamma_\alpha)\,\otimes
e_\mu)\,\eta^{\mu
  \nu}\,\otimes\,\Omega_{g_2,n_2}(e_\nu\,\otimes\,\bigotimes_{\alpha\in 
  S_2}\,\gamma_\alpha) 
\end{align*}
where the sign $\pm$ is the usual one obtained by applying the permutation
induced by $S$ to $(\gamma_1,\gamma_2,\ldots\,\gamma_n)$ taking into
account the grading of $\{\,\gamma_i\,\}$ and where $\{\,e_\alpha\,\}$ is a
homogeneous basis for $\HH$.

\item[iii] Let $\chi\,:\,\M_{g-1,n+2}\,\to\,\M_{g,n}$ be the canonical map
  corresponding to attaching the last two marked points together then
\begin{equation*}
\chi^*\,\Omega_{g,n}(\gamma_1,\gamma_2,\ldots,\gamma_n)\,=\,
\Omega_{g-1,n+2}\,(\gamma_1,\gamma_2,\ldots,\gamma_n, e_\mu,
e_\nu)\,\eta^{\mu\nu}. 
\end{equation*}
\end{description}

Let $\Gamma$ be a stable graph then there is a canonical map $\rho_\Gamma$ 
obtained by composition of the canonical maps 
\[
\prod_{v\in V(\Gamma)}\,\M_{g(v),n(v)}\,\to\,
\M_\Gamma\,\to\,\M_{g,n}.
\]
Since the map $\rho_\Gamma$ can be constructed from morphisms in $(ii)$ and
$(iii)$ above, $\Omega_{g,n}$ satisfies a restriction property of the form
\begin{equation} \label{restriction}
\rho_\Gamma^*\Omega_{g,n}\,=\, P_\Gamma(\,\bigotimes_{v\in V(\Gamma)}
\Omega_{g(v),n(v)}\,)
\end{equation}
where 
\[ P_\Gamma\,:\,\bigotimes_{v\in V(\Gamma)} \R_{g(v),n(v)}\, \to\,
\R_{g,n} \]
is the linear map contracting tensor factors of $\HH$ using the metric $\eta$
induced from successive application of equations $(ii)$ and $(iii)$ above.

The collection $\Omega$ forms a (complete) \cft\ $H_\bullet(\M_{g,n})\, \to\,
\End_{g,n}(\HH)$ via $[c]\,\mapsto\,\int_{[c]}\Omega_{g,n}$. It is this
formulation of \cfts\ which arises most naturally in algebraic geometry.

Notice that the definition of a cohomological field theory is valid even when
enlarges the ground ring from $\nc$ to another ring $\mathcal{K}$. 

From now on, we shall restrict our consideration to genus zero \cfts.

\begin{nota}
Let $V$ be a topological space and let $H^\bullet(V,\nc)$ be given a
homogeneous basis $\mathbf{e}\,:=\,\{\,e_\alpha\,\}_{\alpha\in A}$ and let
$e_0$ denote the identity element. Let
$\sf\,:=\,\{\,s_r^\alpha\,|\,r\geq -1, \alpha\in A\,\}$ be a
collection of formal variables with grading $|s_r^\alpha|
\,=\,2r+|e_\alpha|$. All formal power series and polynomials in a collection
of variables (\emph{e.g.}\ $\sf$) are in the $\nz_2$-graded sense.
\end{nota}

\begin{df}
Let $\Lambda$ consist of formal symbols $q^\beta$ for all
$\beta\,\in\,B(V)$ together with the multiplication
$(q^\beta\,q^{\beta'})\,\mapsto\, q^{\beta+\beta'}$. Let
$\nc[[\Lambda]]$ consist of formal sums
$\sum_{\beta\,\in\,B(V)}\,a_\beta\,q^\beta$ where $a_\beta$ are
elements in $\nc$. Assign to each $q^\beta$, the degree $-2
c_1(V)\,\cap\,\beta$. The product is well-defined according to
\cite[Prop.~II.4.8]{Kol}. This endows $\Lambda$ with the structure of
a semigroup with unit. Furthermore, let $\nc[[\Lambda,\sf]]\, := \,
\nc[[\Lambda]][[\sf]]$, formal power series in the variables $\sf$
with coefficients in $\nc[[\Lambda]]$.
\end{df}

Suppose that $H_2(V,\nz)$ has no torsion then the semigroup $\Lambda$ can be
 realized concretely by choosing a basis  $\{\,e^\epsilon\,\}$ for
$H_2(V,\nz)$ (say the one which is a subset of the basis dual via the
Poincar\'e pairing to $\{\,e_\alpha\,\}$), introducing a formal variable
$q_\epsilon$ for each such basis element, writing
$\beta\,=\,\sum_{\epsilon}\, \beta_\epsilon e^\epsilon$, and defining
$q^\beta\,=\, \prod_{\epsilon}\,q_\epsilon^{\beta_\epsilon}$. The variable
$q_\epsilon$ is then assigned the degree $-2 c_1(V)\,\cap\,e^\epsilon$. 
If $H_2(V,\nz)$ contains a torsion subgroup decomposed into a product of
 cyclic groups then one can do the same by introducing a ``root of unity''
 for each generator of a cyclic subgroup in the obvious manner.

\begin{thm}
Let $V$ be a smooth, projective, convex variety.  For each $n\geq 3$, let
$\Omega_{0,n}$ be elements in $\R_{0,n}(V)[[\Lambda,\sf]]$ defined by
\begin{equation} \label{omega}
\Omega_{0,n}(\gamma_1,\ldots,\gamma_n)\,:=\,\sum_\beta\,\st_*(\,ev_1^*
\gamma_1\,\cdots\, ev_n^*\gamma_n\,\exp(\mathbf{s} \kaf))
\,q^\beta 
\end{equation}
where $\st\,:\,\M_{0,n}(V,\beta)\,\to\,\M_{0,n}$ is the forgetful map,
$\gamma_1, \gamma_2, \ldots, \gamma_n$ are elements in $H^\bullet(V,\nc)$.
The $\Omega\,:=\,\{\,\Omega_{0,n}\,\}$ endows
$(H^\bullet(V,\nc[[\Lambda,\sf]]),\eta)$ with the structure of a \cft\
where $\eta$ is the Poincar\'e pairing extended $\nc[[\Lambda,\sf]]$-linearly.
In particular, if $H^\bullet(V,\nc)$ consists of only even dimensional
cohomology classes then for all values of $\{\,s_r^\alpha\,\}$,
$\Omega\,:=\,\{\,\Omega_{0,n}\,\}$ endows
$(H^\bullet(V,\nc[[\Lambda]]),\eta)$  with the structure of a \cft\ where
$\eta$ is the Poincar\'e pairing.
\end{thm}

\begin{proof}
Let $n\geq 3$, and let $\Gamma$ be a stable graph for $\M_{0,n}$. Let
$\G$ be the set of graphs $\Gamma'$ obtained by decorating the
vertices of $\Gamma$ with elements $\beta(v) \in H_2(V)$ such that
$\sum_{v\in V(\Gamma)} \beta(v) = \beta$. Clearly, each $\Gamma' \in
\G$ determines a stratum in $\M_{0,n}(V,\beta)$. It follows from the
deformation theory considerations in \cite{FP} that the union of
$\M_{\Gamma'}$, $\Gamma' \in \G$, is the scheme theoretic preimage of
$\M_\Gamma$ under the morphism $\st$. 

Consider the following commuting diagram:
\[
\begin{CD}
\coprod_{\Gamma' \in \G} \prod_{v\in V(\Gamma')} \M(v) 
@< j_{\Gamma'} << \coprod_{\Gamma' \in \G} \fprod_{v\in
V(\Gamma')} \M(v) @> \rho_{\Gamma'} >> \M_{0,n}(V,\beta) \\ 
@VV \widetilde{\st} V @VV \widetilde{\st} V @VV \st V \\
\prod_{v\in V(\Gamma)} \M(v) @<\mathrm{id}<<    \prod_{v\in V(\Gamma)}
\M_{0,n(v)}  @>\rho_{\Gamma} >> \M_{0,n}. \\ 
\end{CD}
\]

One has to check that the element $\Omega_{0,n} (\gamma_1, \ldots,
\gamma_n)$ defined by \eqref{omega} possesses the restriction property
\eqref{restriction}. It suffices to check it for each $\beta$. Since
the right square of the commutative diagram is the normalization of a
pull-back square one has that $\rho_\Gamma^*\, \st_* = \sum_{\Gamma' \in
\G}\, \widetilde{\st}_{*} \,j_{\Gamma' *} \,\rho_{\Gamma'}^*$
\cite[6.2]{Fu}. The 
required restriction property follows from the following facts.
Firstly, for each class $\kappa = \kappa_{a,\alpha}$ one has
$\rho_{\Gamma'}^*\, \kappa = \sum_{v\in V(\Gamma')} j_{\Gamma'}^*
\,\kappa(v)$, where $\kappa(v)$ is the corresponding class on $\M(v)$,
as shown in Lem.~\ref{kappa:restr}. Secondly, $\rho_{\Gamma'}^*
ev_i^*\, \gamma_i\, =\, j_{\Gamma'}^* \, ev_{(v),i}^*
(\gamma_i)$, where $v$ is the vertex of $\Gamma'$ to which the tail labeled
$i$ attached, and $ev_{(v),i}$ is the corresponding evaluation
morphism. Thirdly, the composition $j_{\Gamma'*}\, j_{\Gamma'}^*$ is the
multiplication of by the Thom class of the fibered product in the
corresponding direct product. 

Finally, suppose that $H^\bullet(V)$ consists of only even dimensional
classes. By dimensional considerations and the fact that
$\kappa_{0,0}\,=\,n-2$ and $\kappa_{-1}(\gamma)\,=\,\int_\beta\,\gamma$ when
$\gamma$ has degree $2$, the coefficient of $q^\beta$ in $\Omega_{0,n}$
converges when one plugs in numbers for  $s_n^\alpha$ for all
$n\,=\,-1,0,\ldots$ and $\alpha\,=\,0,1,\ldots , r$.
\end{proof}

There is a third way of describing \cfts\ and this is in terms of a
generating function called its potential. We will see in the next section 
that when considering generating functions for the intersection numbers
involving combinations of $\psi$ classes as well, this alternate description
of \cfts\ will be more useful.

\begin{df}
The \emph{potential} $\Phi$ of a genus zero \cft\ is given by a linear map
$\mu\,:\,H_\bullet(\M_{0,\bullet},\nc)\,\to\,\End(\HH)$ of rank $r+1$ is
defined by choosing a basis $\{\,e_0,\ldots,e_r\,\}$ for $\HH$ and let
$\bx\,=\,\sum_{\alpha\,=\,0}^r\, x^\alpha\,e_\alpha$ represent an arbitrary
cohomology class in $\HH$.  Let $I_{0,n}(\bx,\bx,\cdots,\bx)$ be the element
in $\nc[[x^0,x^1\,\ldots,x^r]]$, the ring of formal power series with $\nc$
coefficients in $x^0,\ldots,x^r$ in the graded sense, obtained by using $h$
to pair $\mu_{0,n}([\M_{0,n}])$ with $\bx\otimes \bx \otimes\ldots\otimes
\bx$ and  
\[
\Phi(\bx) := \sum_{n=3}^\infty
\,I_{0,n}(\bx,\bx,\ldots,\bx) \frac{1}{n!}
\]
which is regarded as an element in $\nc[[x^0,\ldots,x^r]]$.
\end{df}

This definition generalizes to arbitrary genera in the obvious manner. The
potential completely characterizes the cohomological field theory in genus
zero as can be seen from the following theorem.

\begin{thm}
\label{thm:wdvvg}
An element $\Phi$ in $\nc[[x^0,\ldots,x^r]]$ is the potential of a
rank $r$, genus zero \cft\ $(\HH,\eta)$ if and only if \cite{KM1,Ma} it
contains only terms which are cubic and higher order in the variables
$x^0,\ldots,x^r$ and it satisfies the WDVV equation
\[
(\partial_{a} \partial_{b} \partial_{e}\Phi)\, \eta^{ef}\, (\partial_{f}
\partial_{c} \partial_{d}\Phi)\, = \,(-1)^{|x_a|(|x_b| + |x_c|)}\,
(\partial_{b} \partial_{c}
\partial_{e} \Phi)\, \eta^{ef} \,(\partial_{f} \partial_{a}
\partial_{d}\Phi),
\]
where $\eta_{ab} := \eta(e_a,e_b)$, $\eta^{ab}$ is in inverse matrix to
$\eta_{ab}$,  $\partial_a$ is derivative with respect to $x^a$, and the
summation convention has been used.
\end{thm} 

The theorem is a consequence of the work of Keel \cite{Ke} who proved that
all relations between boundary divisors in $H_\bullet(\M_{0,n})$ arise from
lifting the basic codimension one relation on $\M_{0,n}$.

As before, one can extend the ground ring $\nc$ above to $\nc[[\Lambda,\sf]]$
in the definition of the potential of a genus zero \cft\  and the above
theorem extends, as well. In our setting, the potential is a formal function
on $\HH\,:=\,H^\bullet(V,\nc[[\Lambda,\sf]])$ of the following kind. $\Phi$
belongs to $\nc[[\Lambda,\sf]][[x^0,\ldots,x^r]]$ where we have used 
\[
I_{0,n}(e_{\alpha_1},\ldots,e_{\alpha_n})\,=\,\int_{\M_{0,n}}\,\left<
\,\Omega_{0,n}\,,\, e_{\alpha_1}\otimes e_{\alpha_2}\otimes\cdots\otimes
e_{\alpha_n}\,\right>.
\]
Again, if $H^\bullet(V)$ consists entirely of even dimensional classes then
plugging in numbers for all $s_n^\alpha$ where $n=-1,0,1,\ldots$ and
$\alpha\,=\,0,1,\ldots,r$, one obtains families of \cft\ structures on
$H^\bullet(V,\nc[[\Lambda]])$.


\section{Topological Recursion Relations}
\label{trr}

In this section we will derive differential equations for the
generating function which incorporates the intersection numbers of the
$\psi$ and $\kappa$ classes. 

\begin{nota}
Let $\SC_{k}^r$ be the set of infinite sequences of non-negative
integers
\[
\mf= (m_{k}^0, m_{k}^1, \ldots, m_{k}^r, m_{k+1}^0, \ldots, 
m_{k+1}^r, m_{k+2}^0, \ldots ) 
\]
such that $m_{a}^\alpha = 0$ for all $a$ sufficiently large. We denote
by $\delf_a^\alpha$ the infinite sequence whose only non-zero entry is
$m_a^\alpha=1$. We use the Latin characters for the lower index and
the Greek characters for the upper index. We will use the notation of
the type
\[
\mf ! := \prod_{\alpha=0}^r \prod_{a\ge k} m_{a}^\alpha ! \quad
\text{and} \quad 
\binom{\mf}{\mf'} := \frac{\mf !}{\mf'\, !\, (\mf - \mf')\, !}
\]
\end{nota}

We also consider $\tf = \{ t_d^\mu \}$, $d\ge 0$, $\mu= 0, \ldots, r$,
$\sf = \{ s_a^\alpha \}$, $a\ge -1$, $\alpha = 0, \ldots, r$, be a
collection of independent formal variables with grading $|t_a^\alpha|\,=\, 2
a - 2 + |e_\alpha|$ and $|s_a^\alpha|\, =\, 2a + |e_\alpha|$.

\subsection{Differential Equations}
We start with introducing the notation for the intersection numbers.
This notation generalizes the one from \cite{W,KK}. 

\begin{df} Let $\beta \in H_2 (V,\nz)$, and $\{ e_\alpha \}$, $\alpha
= 0, \ldots, r$ is a basis of $H^\bullet (V)$. If $\beta \ne 0$, or
$\beta=0$ and $n\ge 3$, then 
\begin{align*}
&\la \tau_{d_{1}}^{\mu_{1}} \tau_{d_{2}}^{\mu_{2}} \ldots
\tau_{d_n}^{\mu_n } \kappa_{a_1, \alpha_1}, \kappa_{a_2, \alpha_2}
\ldots \kappa_{a_l, \alpha_l} \ra_\beta :=  \\ & \hspace{.3in}
\int_{\M_{0,n}(V,\beta)} 
t_{d_{1}}^{\mu_{1}} \psi_{1}^{d_{1}} ev_{1}^*(e_{\mu_{1}}) \,
t_{d_{2}}^{\mu_{2}}\psi_{2}^{d_{2}} ev_{2}^*(e_{\mu_{2}}) \ldots
t_{d_n}^{\mu_n}\psi_n^{d_n} ev_n^*(e_{\mu_n}) \times \\ & 
\hspace{.8in} 
s_{a_1}^{\alpha_1} \kappa_{a_1, \alpha_1} \,
s_{a_2}^{\alpha_2} \kappa_{a_2, \alpha_2} \ldots 
s_{a_l}^{\alpha_l}\kappa_{a_l,\alpha_l} =
(-1)^\epsilon  t_{d_{1}}^{\mu_{1}} \ldots t_{d_n}^{\mu_n}\,
s_{a_1}^{\alpha_1} \ldots s_{a_l}^{\alpha_l} \times
\\ & \hspace{.3in}
\int_{\M_{0,n}(V,\beta)} 
\psi_{1}^{d_{1}} ev_{1}^*(e_{\mu_{1}}) \ldots 
\psi_n^{d_n} ev_n^*(e_{\mu_n}) \,
\kappa_{a_1, \alpha_1} \ldots \kappa_{a_l,\alpha_l}. 
\end{align*}
If $\beta=0$ and $n < 3$ we define all intersection numbers to be
equal to $0$. 

Let the sequence $(d_{1},\mu_{1}), (d_{2}, \mu_{2}), \ldots, (d_n,
\mu_n)$ contain $m_{d}^\mu$ pairs $(d,\mu)$, where $d\ge 0$, $0\le \mu
\le r$, and the sequence $(a_1, \alpha_1), (a_2, \alpha_2), \ldots,
(a_l, \alpha_l)$ contain $p_{a}^{\alpha}$ pairs $(a,\alpha)$ where
$a\ge -1$, $0\le \alpha \le r$. Then we also denote the intersection
number above by $\la \tauf^\mf \kaf^\pf \ra_\beta$.

We also set 
\[
\la \tauf^\mf \kaf^\pf \ra := 
\sum_\beta \la \tauf^\mf \kaf^\pf \ra_\beta \, q^\beta, 
\]
where $q$ is the formal variable introduced in the previous section.
\end{df}

\begin{rem}
We put the formal variables in the definition of the intersection
numbers in order to take care of the signs provided $V$ has odd
cohomology classes. If $V$ has only even cohomology classes, then
$\epsilon$ in the formula above is $0$. Then setting all formal $t$
variables equal to $1$ and all formal $s$ variables equal to $0$ we
obtain the definition due to Witten \cite{W}. In this paper we assume
that $V$ may have odd cohomology classes even though we do not know
examples of convex manifolds with odd cohomology classes.
\end{rem}

\begin{df}
We define $H(\tf; \sf)$ in $\nc[[\Lambda,\tf,\sf]]$ by
\[
H(\tf; \sf) := 
\sum_{\mf \in \SC_{0}^r\ \pf \in \SC_{-1}^r}
\frac{1}{\mf !}\, \frac{1}{\pf !}\, 
\la \tauf^\mf \kaf^\pf \ra 
\]
\end{df}
Notice that $H$ has grading $2(\dim_\nc V -3)$.

We are going to exploit the presentations of the $\psi$ classes and
the $\kappa$ classes in Lem.~\ref{psi:pres} and Cor.~\ref{kappa:pres}
in order to obtain a system of differential equations for $H$. The
explicit presentation of the tautological classes allows one to obtain
recursion relations for the intersection numbers, and these recursion
relations in  turn provide the corresponding differential
equations. The technical details are standard (cf.~\cite{W,KMZ,KK}),
and we omit them.

We recall that we denote by $\eta$ the Poincar\'e pairing on
$H^\bullet (V)$. We denote $\eta (e_{\alpha_1}, e_{\alpha_2})$ by
$\eta_{\alpha_1, \alpha_2}$, and by $\eta^{\alpha_1, \alpha_2}$ the
coefficients of the inverse metric. We denote by
$c^{\alpha_1}_{\alpha_2, \alpha_3}$ the coefficients of the
multiplication tensor, that is
\[
e_{\alpha_2} e_{\alpha_3} = 
\sum_{\alpha_1=0}^r c^{\alpha_1}_{\alpha_2, \alpha_3} e_{\alpha_1}.
\]

\begin{prop}
\label{prop:trr}
The intersection numbers defined above satisfy the recursion relations
below called the topological recursion relations. The first two
relations are valid when $a_1 \ge 1$.
\begin{align*}
& (t_{a_1}^{\alpha_1} )^{-1} \la \tauf^{\mf + \delf_{a_1}^{\alpha_1} + 
\delf_{a_2}^{\alpha_2} + \delf_{a_3}^{\alpha_3} }   
\kaf^\pf  \ra_\beta = 
(t_{a_1-1}^{\alpha_1} )^{-1} (t_{0}^{\sigma_1} )^{-1} 
(t_{0}^{\sigma_2} )^{-1} \times
\\& \hspace{.2in}
\sum\binom{\mf}{\mf'} \binom{\pf}{\pf'} 
\la \tauf^{\mf'+ \delf_{a_{1} -1 }^{\alpha_1} + 
\delf_{0}^{\sigma_1} } 
\kaf^{\pf'} \ra_{\beta_1}
\eta^{\sigma_1,\sigma_2} 
\la \tauf^{\mf'' + \delf_{0}^{\sigma_2} + 
\delf_{a_2}^{\alpha_2} + \delf_{a_3}^{\alpha_3} } 
\kaf^{\pf''} \ra_{\beta_2}  \\
& (s_{a_1}^{\alpha_1} )^{-1} \la \tauf^{\mf + 
\delf_{a_2}^{\alpha_2} + \delf_{a_3}^{\alpha_3} }   
\kaf^{\pf + \delf_{a_1}^{\alpha_1} }  \ra_\beta = 
(s_{a_1-1}^{\alpha_1} )^{-1} (t_{0}^{\sigma_1} )^{-1} 
(t_{0}^{\sigma_2} )^{-1} \times
\\& \hspace{.2in}
\sum \binom{\mf}{\mf'} \binom{\pf}{\pf'} 
\la \tauf^{\mf'+ \delf_{0}^{\sigma_1} } 
\kaf^{\pf' + \delf_{a_{1} -1 }^{\alpha_1} } \ra_{\beta_1}
\eta^{\sigma_1,\sigma_2} 
\la \tauf^{\mf'' + \delf_{0}^{\sigma_2} + 
\delf_{a_2}^{\alpha_2} + \delf_{a_3}^{\alpha_3} } 
\kaf^{\pf''} \ra_{\beta_2}  \\
& (s_{0}^{\alpha_1} )^{-1} \la \tauf^{\mf + 
\delf_{a_2}^{\alpha_2} + \delf_{a_3}^{\alpha_3} }   
\kaf^{\pf + \delf_{0}^{\alpha_1} }  \ra_\beta = 
(s_{-1}^{\alpha_1} )^{-1} (t_{0}^{\sigma_1} )^{-1} 
(t_{0}^{\sigma_2} )^{-1} \times 
\\& \hspace{.2in}
\sum \binom{\mf}{\mf'} \binom{\pf}{\pf'} 
\la \tauf^{\mf'+ \delf_{0}^{\sigma_1} } 
\kaf^{\pf' + \delf_{ -1 }^{\alpha_1} } \ra_{\beta_1}
\eta^{\sigma_1,\sigma_2} 
\la \tauf^{\mf'' + \delf_{0}^{\sigma_2} + 
\delf_{a_2}^{\alpha_2} + \delf_{a_3}^{\alpha_3} } 
\kaf^{\pf''} \ra_{\beta_2} + \\ & \hspace{.2in}
\sum_{a=0}^{\infty}
\sum_{\alpha, \nu=0}^r m_a^\alpha c^\nu_{\alpha, \alpha_1} 
t_a^\alpha (t_a^\nu)^{-1} 
\la \tauf^{\mf - \delf_{a}^{\alpha} + \delf_{a}^{\nu}  + 
\delf_{a_2}^{\alpha_2} + \delf_{a_3}^{\alpha_3} } 
\kaf^\pf  \ra_\beta
\end{align*}
where the summation is over $\beta_1 + \beta_2 = \beta,
\mf'+\mf''=\mf, \pf'+\pf''=\pf$, and $\sigma_1, \sigma_2 $ varying
from $0$ to $r$.

Equivalently, the function $H(\tf;\sf)$ satisfies the system of
differential equations below. The first two equations are valid when
$a_1 \ge 1$.
\begin{align*}
\frac{\dd^3 H}{\dd t_{a_1}^{\alpha_1} \dd t_{a_2}^{\alpha_2} 
\dd t_{a_3}^{\alpha_3} }  = 
&\sum_{\sigma_1, \sigma_2} 
\frac{\dd^2 H}{\dd t_{a_1-1}^{\alpha_1} \dd t_{0}^{\sigma_1} }
\, \eta^{\sigma_1,\sigma_2} \,
\frac{\dd^3 H}{\dd t_{0}^{\sigma_2} \dd t_{a_2}^{\alpha_2} 
\dd t_{a_3}^{\alpha_3} } \\ 
\frac{\dd^3 H}{\dd s_{a_1}^{\alpha_1} \dd t_{a_2}^{\alpha_2} 
\dd t_{a_3}^{\alpha_3} }  =
&\sum_{\sigma_1, \sigma_2} 
\frac{\dd^2 H}{\dd s_{a_1-1}^{\alpha_1} \dd t_{0}^{\sigma_1} }
\, \eta^{\sigma_1,\sigma_2} \,
\frac{\dd^3 H}{\dd t_{0}^{\sigma_2} \dd t_{a_2}^{\alpha_2} 
\dd t_{a_3}^{\alpha_3} } \\ 
\frac{\dd^3 H}{\dd s_{0}^{\alpha_1} \dd t_{a_2}^{\alpha_2} 
\dd t_{a_3}^{\alpha_3} }  = 
&\sum_{\sigma_1, \sigma_2} 
\frac{\dd^2 H}{\dd s_{-1}^{\alpha_1} \dd t_{0}^{\sigma_1} }
\, \eta^{\sigma_1,\sigma_2} \,
\frac{\dd^3 H}{\dd t_{0}^{\sigma_2} \dd t_{a_2}^{\alpha_2} 
\dd t_{a_3}^{\alpha_3} } \\+
&\sum_{a=0}^\infty \sum_{\alpha, \nu=0}^r 
c^{\nu}_{\alpha, \alpha_1} t_{a}^{\alpha} 
\frac{\dd^3 H}{\dd t_{a}^{\nu} \dd t_{a_2}^{\alpha_2} 
\dd t_{a_3}^{\alpha_3} },
\end{align*}
where the indices $\sigma_1, \sigma_2$ vary from $0$ to $r$. \qed
\end{prop}

\begin{rem}
If one considers the function $F(\tf) := H(\tf;\zef)$, then the first
equations reduce to the equations first written by Witten in \cite{W}.
Witten showed \cite{W} that the WDVV equations in Th.~\ref{thm:wdvvg}
follow from the differential equation for $F(\tf)$.
\end{rem}

\begin{crl}
Let $H$ be as above. Let $x^\alpha\,:=\,t_0^\alpha$,
$\bx\,:=\,(x^0,\ldots,x^r)$, and let $\tf\,:=\,(\bx,\widetilde{\tf})$
where $\widetilde{t}$ consists of all $t_n^\alpha$ such that $n\,\geq\,1$. Let
$\Phi(\bx)\,:=\,H(\bx,\widetilde{\tf};\sf)$ be regarded as an element of
$\nc[[\Lambda,\widetilde{\tf},\sf]][[\bx]]$ then $\Phi$ endows
$H^\bullet(V,\nc[[\Lambda, \widetilde{\tf}, \sf]])$ with the structure of a
\cft\  (over the ground ring $\nc[[\Lambda,\widetilde{\tf},\sf]]$). 
\end{crl}
\begin{proof}
Our proof follows that in \cite{W} but keeping track of signs. Differentiate
the first set of equations in the last theorem with respect to the variable
$t_{a_4}^{\alpha_4}$, setting $a_1\,=1\,$ and $a_2\,=\,a_3\,=\,a_4\,=0$ and
noting that the left hand side is essentially (graded) symmetric under
permutation of the subscripts $2,3,4$. Therefore, the antisymmetric part of
the right hand side must vanish. This is nothing more than the WDVV equations
(with the correct signs).  This implies (by Thm \ref{thm:wdvvg}) that $\Phi$
is a potential for a genus zero \cft. 
\end{proof}

\subsection{Puncture and Dilaton Equations} Here we introduce
analogues of the puncture and the dilaton equations. Using these
equations and the topological recursion relations we will show that
the intersection numbers defined above can be expressed in terms of
Gromov--Witten invariants.

For the rest of this subsection let us fix the following notation. Let
$\pi: \M_{0,n+1} (V,\beta) \to \M_{0,n} (V, \beta)$ denote the
universal stable map which ``forgets'' the $n+1^{\text{st}}$ marked
point. To distinguish the $\psi$ and the $\kappa$ classes upstairs and
downstairs we put ``hats'' over the former. Since the pull back of
$ev_i (e_\alpha)$ is still $ev_i (e_\alpha)$ we denote these classes
by the same symbol unstairs and downstairs. 

In addition, if $\mf = \{ m_a^\alpha \} \in \SC_{k}^r$, $k=-1,0$, then
\begin{equation*}
|\mf| := \sum_{\alpha=0}^r \sum_{a\ge 0} a m_a^\alpha, 
\quad  \|\mf\| := \sum_{\alpha=0}^r \sum_{a\ge 0} m_a^\alpha
\end{equation*}
and if $\mf = \{ m_a^\alpha \} \in \SC_{0}^r$, then 
\begin{align*}
e(\mf) &:= e_0^{m_0^0} \ldots e_r^{m_0^r} \,
e_0^{m_1^0} \ldots e_r^{m_1^r}\, e_0^{m_2^0} \ldots \\
\sf^\mf &: = \ldots (s_2^0)^{m_2^0} \,
(s_1^r)^{m_1^r} \ldots (s_1^0)^{m_1^0} \,
(s_0^r)^{m_0^r} \ldots (s_0^0)^{m_0^0}
\end{align*}
Note that in both cases we have summation over $a\ge 0$, not over
$a\ge k$, and that the order of the terms in the products is the
opposite. We denote by $c_{e(\mf),\mu}^\nu$ the coefficients
of $e_\nu$ in the product $e (\mf) e_\mu$. We also say that
$\mf \in \SC^r_{0,0}$ if $\mf \in \SC_0^r$, and $m_a^\alpha = 0$
provided $a\ge 1$.

It is easy to derive the following formulas for the pull back of the
$\psi$ and $\kappa$ classes with respect to $\pi$ (cf.~\cite{Pa}).
Namely,
\[
\psio_i^{a} = \pi^* \psi_i^a + \pi^* \psi_i^{a-1} D_i \quad 
\text{and} \quad
\kao_{a,\alpha} = \pi^* \kappa_{a,\alpha} + 
\psio_{n+1}^a ev_{n+1}^* (e_\alpha). 
\]

The pull back formulas above induce the recursion relations
\eqref{dilaton} below. These are analogues of the puncture and dilaton
equations. The proof is standard, and we omit the details
(cf.~\cite[Sec.~3.1]{KK}). The variables take care of the signs as in
the case of the topological recursion relations. When $a\ge 1$, or
when $a=0$ and $\mf \in \SC_{0,0}^r$, one has
\begin{multline}
\label{dilaton}
(t_a^\alpha)^{-1} 
\la \tauf^{\mf + \delf_a^\alpha} \kaf^\pf \ra_\beta = \\  \hspace{0.2in}
\sum_{\substack{\pf' + \pf'' = \pf \\ \pf'' \in \SC_0^r} }
\sum_{\nu=o}^r
\binom{\pf}{\pf'} c_{ e(\pf''),\alpha}^\nu 
\sf^{\pf''} (s_{|\pf''| + a -1}^\nu)^{-1} 
\la \tauf^\mf \kaf^{\pf' + \delf^\nu_{|\pf''| + a -1}} \ra_\beta. 
\end{multline}
Here we take $\pf'' \in \SC_{0,r}$ since $\pi^* (\kappa_{-1,\alpha}) =
\kao_{-1,\alpha}$. Of course, one can not apply \eqref{dilaton} when
$\beta=0$ and $\| \mf + \delf_a^\alpha \| = 3$ since the space
$\M_{0,2}(V,0)$ does not exist. One can also derive a recursion
relation in case when $a=0$, and $\mf$ is arbitrary, but it is rather
messy, and we will not need it in the sequel. The recursion relations
\eqref{dilaton} are equivalent to the differential equations
\begin{equation*}
\frac{\dd (H - H_{in}) }{\dd t_{a}^{\alpha} } =
\sum_{\pf \in \SC_0^r} \sum_{\nu=0}^r 
c_{ e(\pf),\alpha}^\nu 
\frac{ \sf^\pf}{\pf !} \,
\frac{\dd H} {\dd s_{| \pf | + a -1 }^{\nu}} 
\end{equation*}
valid for all $a\ge 1$, and for $a=0$ if instead of $H$ we consider
the function $H_1$ which is obtained from $H$ by setting all variables
$t_i^\mu$ with $i\ge 1$ to zero. The function $H_{in}$ is the part of
$H$ which corresponds intersection numbers on $\M_{0,3}(V,0)$. It is a
function in variables $t_0^\alpha$, $s_0^\alpha$, $\alpha=0, \ldots,
r$, determined by the cup-product on $V$ since $\kappa_{0,\alpha}$ on
$\M_{0,3}(V,0)$ is equal to $ev^* (d_\alpha)$. (All evaluation maps
are equal when $\beta=0$.) 

We want to show that \eqref{dilaton} and Prop.~\ref{prop:trr} allow us
to compute inductively all intersection numbers of the $\psi$ and the
$\kappa$ classes provided we know the Gromov--Witten invariants of
$V$. Notice that \eqref{dilaton} allows one to eliminate
$\tau_a^\alpha$ with $a\ge 1$ from the intersection numbers. Therefore
we reduced the calculation to the case of the intersection numbers of
the $\kappa$ classes and the pull backs of the cohomology classes from
$V$ if no $\psi$ classes are present.

First we want to reduce everything to the case of the intersection
numbers of $\kappa_{-1,\alpha}$ and $ev_i^* (e_\mu)$. Suppose that we
have an intersection number of the type $\la \tauf^\mf \kaf^\pf
\ra_\beta$, where $\mf \in \SC_{0,0}^r$, and there exists $p_a^\alpha
>0$ for some $a \ge 0$. If $\| \mf \| \ge 2$, then the second or third
recursion relation in Prop.~\ref{prop:trr} allows us to express it in
terms of the intersection numbers of the type $\la \tauf^{\mf'}
\kaf^{\pf'} \ra_{\beta_1}$, where $\mf' \in \SC_{0,0}^r$, $\| \mf' \|
\ge 1$, and either $\| \pf' \| < \| \pf \|$, or $\| \pf' \| = \| \pf
\|$, and $| \pf' | < | \pf |$.

However, if $\| \mf \| = 1$, we can not apply the recursion relations
from Prop.~\ref{prop:trr}. Here we have to use \eqref{dilaton}.
Namely, pick $\alpha$ such that $| e_\alpha |=2$ and $\int_\beta \alpha
\ne 0$. (Clearly, such $\alpha$ exists if $V$ is not a point.) When we
apply \eqref{dilaton} to $\la \tauf^{\mf + \delf_0^\alpha} \kaf^\pf
\ra_\beta$ we can split the sum in the right hand side into two parts.
One term corresponds to $\pf' = \pf$, and $\pf'' = \zef$. This term is
\[
\la \tauf^\mf \kaf^\pf \kappa_{-1,\alpha} \ra_\beta\, = \,
 \la \tauf^\mf \kaf^\pf\ra_\beta\, \int_\beta \alpha \,
\]
Here and in the rest of the argument we disregard the formal variables
since they just take care of the signs which do not play an important role in
the our argument. The second term is 
\[
\sum_{\substack{\pf' + \pf'' = \pf \\ \zef \ne \pf'' \in \SC_0^r} }
\sum_{\nu=0}^r
\binom{\pf}{\pf'} c_{ e(\pf''),\alpha}^\nu 
\la \tauf^\mf \kaf^{\pf' + \delf^\nu_{|\pf''|-1}} \ra_\beta.
\]
For all these terms  
\[
\| \pf'  + \delf^\nu_{|\pf''|-1} \| =
\| \pf' \| + 1 \le  \| \pf \|. 
\]
The equality is attained only when $\pf'' = \delf_{j}^\mu$ for $j\ge
1$. But for such terms $| \pf' + \delf^\nu_{|\pf''|-1} | < | \pf |$.
It follows that we can express each intersection
number $\la \tauf^\mf \kaf^\pf \ra_\beta$ with $\mf \in \SC_{0,0}^r$,
$\| \mf \| \ge 1$ in terms of the intersection numbers $\la
\tauf^{\mf'} \kaf^{\pf'} \ra_{\beta_1}$, where either $\| \pf' \| < \|
\pf \|$, or $\| \pf' \| = \| \pf \|$ and $| \pf' | < | \pf |$. Note,
that it is easy to extend our argument to the case when $\| \mf \| =
0$. 

Inductively this reduces the computation of the intersection numbers
to the integrals of the type
\begin{equation}
\label{growit}
\int_{\M_{0,n}(V,\beta)} ev_1^* (\gamma_1) \ldots ev_n^* (\gamma_n)
\kappa_{-1}(\gamma_{n+1}) \ldots \kappa_{-1} (\gamma_{n+l}),
\end{equation}
where $\gamma_i$ are cohomology classes of $V$. Recall that $\pi^*
(ev_i^* (\gamma)) = ev_i^*(\gamma)$ and $\pi^* (\kappa_{-1}(\gamma)) =
\kappa_{-1} (\gamma)$. Iterating this equality one gets that the
intersection number \eqref{growit} is equal to
\begin{equation*}
\int_{\M_{0,n+l}(V,\beta)} ev_1^* (\gamma_1) \ldots ev_n^* (\gamma_n)
ev_{n+1}^* (\gamma_{n+1}) \ldots ev_{n+l}^* (\gamma_{n+l}). 
\end{equation*}
The integral above is the Gromov--Witten invariant evaluated on the
classes $\gamma_1, \ldots, \gamma_{n+l}$. 

\begin{rem}
In \cite{KM2} the authors show how to compute the intersection numbers
of the $\psi$ classes and the pull backs of the cohomology classes
from $V$ provided that the Gromov--Witten invariants are known. One
can show that \eqref{dilaton} allows one to eliminate the $\kappa$ classes
from the intersection numbers. This provides an alternative way of
calculating these intersection numbers via the Gromov--Witten
invariants. 
\end{rem}

\subsection{$\CP{1}$: An Example}

The simplest application of the above equations is to the case where $V$ is a
point.  This situation is discussed in some detail in \cite{KK} in genus zero
and one. 

The next simplest case is when $V\,=\,\CP{1}$. Let $\{\,e_0,\,e_1\,\}$ be a
basis for $H^\bullet(\CP{1})$ where $e_0$ is the identity element and $e_1$
is the element in $H^2(\CP{1})$ such that the Poincar\'e metric
$\eta_{0,1}\,:=\,1$. By abuse of notation, an element $\beta$ in
$H_2(\CP{1},\nz)$ is sometimes regarded as an integer (also denoted by
$\beta$) times the fundamental class of $\CP{1}$. The moduli space of stable
maps $\M_{0,n}(V,\beta)$, has complex dimension $-2 + n + 2\beta$. 
On $\M_{0,n}(V,\beta)$, there are the identities $\kappa_{0,0}\,=\,(n-2) e_0$
and $\kappa_{-1,0}\,=\,0$, and $\kappa_{-1,1}\,=\,\beta$.

Let $x^\alpha \, := \, t_0^\alpha$ for $\alpha\,=\,0,1$.  The generating
function $H(\sf,\tf)$ has degree $-4$. We will now present an explicit
formula for $H(x^0,x^1,s_{-1}^1, s_0^0, s_0^1; q)$ where the $q$ dependence
has been made explicit and it is understood that all other variables have
been set to zero. In particular, $H$ contains no $\psi$ classes at all.

The generating function $\Hp$ consists of all terms in $H$ except those
arising from intersection numbers on $\M_{0,3}(V,0)$. One can directly verify
that 
\begin{multline}
H(x^0,x^1,s_{-1}^1, s_0^0, s_0^1; q)\,=\,\Hp(x^0,x^1,s_{-1}^1, s_0^0, s_0^1;
q) \,\\  + \, e^{s_0^0}\,(\,\frac{1}{2}\,(x^0)^2\,x^1 \,+\,
\frac{1}{6}\,(x^0)^3\,s_0^1\,) 
\end{multline}
where the usual cup product on $H^\bullet(\CP{1})$ corresponds to the special
case where all of the $s_a^\alpha$ variables have been set to zero.

The puncture/dilaton equations for $\CP{1}$ acting upon the generating
function $\Hp\,:=\,\Hp(x^0, x^1, s_{-1}^1, s_0^0, s_0^1; q)$ are
\[
\frac{\partial \Hp}{\partial x^0}(x^0, x^1, s_{-1}^1, s_0^0, s_0^1; q) \, =\,
e^{s_0^0}\,s_{0}^1\,\frac{\partial \Hp}{\partial s_{-1}^1}
\]
and
\[
\frac{\partial \Hp}{\partial x^1}(x^0, x^1, s_{-1}^1, s_0^0, s_0^1; q) \, =\,
e^{s_0^0}\,\,\frac{\partial \Hp}{\partial s_{-1}^1}.
\]
The relevant topological recursion relations for the generating function
$H\,:=\,H(x^0,x^1,s_{-1}^1,s_0^0,s_0^1;q)$ becomes
\[
\frac{\partial H_{\alpha_1 \alpha_2}}{\partial s_0^0} \,=\, x^0 H_{0 \alpha_1
	\alpha_2}\,+\, x^1 H_{1 \alpha_1 \alpha_2}
\]
and
\[
\frac{\partial H_{\alpha_1 \alpha_2}}{\partial s_0^1}\,=\, 
\frac{\partial  H_0}{\partial s_{-1}^1}\,H_{1 \alpha_1 \alpha_2}\,+\,
\frac{\partial H_1}{\partial s_{-1}^1}\,H_{0 \alpha_1 \alpha_2}\,+
x^0 H_{1 \alpha_1 \alpha_2}
\]
for all $\alpha_1,\alpha_2\,=\,0,1$. Furthermore, $H_i$ denotes the
derivative of $H$ with respect to $x^i$, $H_{i j}$ is the derivative of $H$
with respect to the variables $x^i$ and $x^j$, and similarly for $H_{i j
  k}$. 

The puncture/dilaton equations yields a solution of the form
\begin{equation}
\label{eq:form}
\Hp(x^0,x^1,s_{-1}^1, s_0^0, s_0^1; q)\,=\,
\sum_{n=1}^\infty\, e^{-2 s_0^0}\,\qh^n\,\frac{(s_0^1)^{2n-2}}{(2n -
  2)!}\,h_n
\end{equation}
where $\qh\,:=\,q \exp(s_{-1}^1\,+\,e^{s_0^0}\,(x^1 + s_0^1 x^0))$ and $h_n$
are numbers. The form of $\Hp(x^0,x^1,s_{-1}^1, s_0^0, s_0^1; q)$ above
implies that it is enough to know the function $\Hp(x^0,x^1,s_0^1;q)$ (where
$s_{-1}^1$ and $s_0^0$ have both been set to zero) since one can recover
$\Hp(x^0,x^1,s_{-1}^1, s_0^0, s_0^1; q)$ from it. \emph{Henceforth, unless
otherwise stated, assume that $H$ and $\Hp$ depend only upon the
variables $x^0,x^1,s_0^1$ and $q$.}

The topological recursion relations can then all be reduced to the single
equation 
\[
\label{eq:penult}
\frac{\partial\,\Hp''}{\partial
  s_0^1}\,=\,2\,s_0^1\,\Hp''\,\Hp'''\,+\,x^0\,\Hp''' 
\]
where $'$ represents the derivative $q\frac{\partial}{\partial
  q}$. 

This equation gives rise to the recursion relations
\[
h_{n+1}\,=\,\sum_{l=1}^n\,\frac{(2n-1)!\, l^2\, (n+1-l)^2}{(2l-2)!\, ( 2 (n-l) )!\, (n+1)}\, h_l\,h_{n+1-l}
\]
and 
$h_1\,=\,1$.

The first few terms of the potential are readily seen to be
\begin{multline}
H(x^0,x^1,s_{-1}^1,s_0^0,s_0^1;q)\,=\,
  \exp(s_0^0)\,(\,\frac{1}{2}\,(x^0)^2\,x^1 + \frac{1}{6}\,(x^0)^3\,s_0^1\,)\\
  +\exp(-2 s_0^0)\,(\qh + \qh^2\,\frac{(s_0^1)^2}{4} +
  \qh^3\,\frac{(s_0^1)^4}{6} + \qh^4\, \frac{(s_0^1)^6}{6} + \mathrm{o}(q^5
  (s_0^1)^8)). 
\end{multline}

As expected, after setting all $s_a^\alpha$ variables to zero one recovers
\[
H(x^0,x^1;q)\,=\,\frac{1}{2}\,(x^0)^2\,x^1 \,+\,q \exp(x^1),
\]
the usual Gromov-Witten potential of $\CP{1}$. More generally,
the potential $H$ parametrizes Frobenius algebras structures on
$H^\bullet(\CP{1})$ which are further deformations of the quantum cup product
on $H^\bullet(\CP{1})$.

It is amusing to note that integrating equation \ref{eq:penult} twice, one
can obtain an expression for $\Hp'$ in the form of a transcendental equation,
namely,  
\[
(1-u)\,\exp(u)\,=\,q (s_0^1)^2\,\exp(1 + s_0^1 x^0 + s^1)
\]
where 
\[
u\,:=\,\sqrt{1 - 2 (s_0^1)^2\,\Hp'(x^0,x^1,s_0^1;q)}.
\]

\bibliographystyle{amsplain}

\providecommand{\bysame}{\leavevmode\hbox to3em{\hrulefill}\thinspace}

\end{document}